\input ./style/arxiv-general.cfg
\documentclass[aop,MSNbibl,dvips]{arximspdf}
\makeatletter
   \@ifpackageloaded{graphicx}{}{\usepackage{graphicx}}
\makeatother

\doi{10.1214/14-AOP975}
\volume{44}
\issue{1}
\pubyear{2016}
\firstpage{324}
\lastpage{359}
\docsubty{FLA}

\makeatletter
\newcommand{\rrvert}{\vert}
\newcommand{\rrVert}{\Vert}
\newcommand{\llvert}{\vert}
\newcommand{\llVert}{\Vert}
\renewcommand{\mid}{|}
\def\var{\operatorname{Var}}
\def\cov{\operatorname{Cov}}
\newtheorem{theorem}{Theorem}
\newtheorem{lemma}[theorem]{Lemma}
\newtheorem{proposition}[theorem]{Proposition}
\newtheorem{corollary}[theorem]{Corollary}
\newproclaim{remark}[theorem]{Remark}
\newcommand{\ud}{\mathrm{d}}
\newcommand{\anne}{\mathbb E_\lambda}
\newcommand{\annp}{\mathbb P_\lambda}
\newcommand{\anp}{\bar{\mathbb{P}}_\lambda}
\newcommand{\ane}{\bar{\mathbb{E}}_\lambda}
\newcommand{\mb}[1]{\mathbb{#1}}
\newcommand{\mc}[1]{\mathcal{#1}}
\newcommand{\bmb}[1]{\bar{\mathbb #1}}
\makeatother

\begin{document}
\begin{frontmatter}

\title{Einstein relation for random walks in
random~environment\thanksref{T1}}
\runtitle{Einstein relation for RWRE}

\begin{aug}
\author[A]{\fnms{Xiaoqin}~\snm{Guo}\corref{}\ead[label=e1]{guoxx097@umn.edu}}
\runauthor{X. Guo}
\affiliation{Technische Universit\"{a}t M\"{u}nchen}
\address[A]{Fakult\"{a}t f\"{u}r Mathematik\\
Technische Universit\"{a}t M\"{u}nchen\\
Boltzmannstr.~3\\
85748 Garching\\
Germany\\
\printead{e1}}
\end{aug}
\thankstext{T1}{Supported by NSF Grant DMS-08-04133 and ERC StG Grant 239990.}

%
\received{\smonth{10} \syear{2013}}
%
\revised{\smonth{7} \syear{2014}}

%
\begin{abstract}
In this article, we consider the speed of the random walks in a
(uniformly elliptic and i.i.d.)
random environment (RWRE) under perturbation. We obtain the derivative
of the speed of the RWRE w.r.t. the perturbation, under the assumption
that one of the following holds: (i) the environment is balanced and
the perturbation satisfies a Kalikow-type ballisticity condition, (ii)
the environment satisfies Sznitman's ballisticity condition. This is a
generalized version of the Einstein relation for RWRE.

Our argument is based on a modification of Lebowitz--Rost's argument
developed in [\textit{Stochastic Process. Appl.} \textbf{54} (1994) 183--196]
and a new regeneration structure for the
perturbed balanced environment.
\end{abstract}

%
\end{frontmatter}

\setcounter{footnote}{1}
\section{Introduction}
In the 1905, Einstein (\cite{Einstein}, pages 1--18) investigated the
movement of suspended particles in a liquid under the influence of an
external force. He established the following mobility--diffusivity relation:
{\renewcommand{\theequation}{ER}
\begin{equation}\label{eqER}
\lim_{\lambda\to0}\frac{v_\lambda}{\lambda}\sim D,
\end{equation}}%
where $\lambda$ is the size of the perturbation, $D$ is the diffusion
constant of the equilibrium state and $v_\lambda$ is the effective
speed of the random motion in the perturbed media.
General derivations of this principle assume reversibility.

Recently, there has been much interest in studying the Einstein
relation for reversible motions
in a perturbed random media, where the perturbation is proportional to
the original environment; see \cite{Le,KO,GMP,BHOZ}. However, it is
not clear whether (\ref{eqER}) still holds in nonreversible set-up, for
example, random walks in random environments (RWRE), and several
interesting questions are either open or not discussed:
is $v_\lambda$ monotone (in an appropriate sense) and differentiable
with respect to $\lambda$?
What if the perturbation of the environment is not propositional to the
original one? If the original environment is ballistic, (\ref{eqER}) is not
expected to hold, but what can we say about the derivative of the velocity?

Motivated by these questions, we study the speed of RWRE under general
perturbations, where the original environment is either balanced or ballistic.
In the balanced case, when the perturbation is proportional to the
original environment, we obtain the Einstein relation. (This result
forms part of the author's doctoral thesis \cite{Guo}.)
Moreover, we provide a new interpretation of this relation. Namely, in
our context, the Einstein
relation is a consequence of the weak convergence of the invariant
measures for the ``environment viewed from the point of view of the
particle'' process, which holds even for more general perturbations that
satisfies a Kalikow-type condition.
In the ballistic case, we can quantify the rate of the weak convergence.
As a corollary, we obtain the derivative of the speed w.r.t. the size
of the perturbation (for both the balanced and the ballistic cases).


We define the model as follows.

An (uniformly elliptic) \textit{environment} $\omega\dvtx \mathbb
{Z}^d\times\{e\in\mathbb{Z}^d\dvtx \llvert e\rrvert=1\}\to
[\kappa, 1)$ is a function that satisfies
\[
\sum_{e\dvtx \llvert e\rrvert=1}\omega(x,e)=1\qquad\forall x\in
\mathbb{Z}^d,
\]
where $\kappa>0$ and $\llvert\cdot\rrvert$ is the $l^2$-norm.
The random walks in the environment $\omega$ starting from $x$ is the
Markov chain $(X_n)_{n\ge0}$ with transition probability $P_\omega$
specified by
\begin{eqnarray*}
P_\omega^x(X_0=x)&=&1,
\\
P_\omega^x(X_{n+1}=y+e\mid X_n=y)&=&
\omega(y,e).
\end{eqnarray*}

Following Sabot \cite{Sabot}, we consider a perturbed environment
\[
\omega^\lambda:=\omega+\lambda\xi, \qquad\lambda\in[0,\kappa/2),
\]
where $\xi\dvtx \mb{Z}^d\times\{e\in\mb{Z}^d\dvtx  \llvert e\rrvert
=1\}\to[-1,1]$ satisfies
\[
\sum_{e\dvtx \llvert e\rrvert=1}\xi(x,e)=0\qquad\forall x.
\]
We denote the local environment at $x$ as $\omega_x:=(\omega
(x,e))_{e\dvtx \llvert e\rrvert=1}$ and write
\[
\zeta:=(\omega,\xi).
\]
We endow the set $\Omega$ of all $\zeta$ with a probability measure
$\mc P$ such that $(\zeta_x)_{x\in\mathbb Z^d}$
are independent and identically distributed ({i.i.d.}).

The measure $P_{\omega^\lambda}^x$ for a fixed $\omega$ is called
the \textit{quenched law}. The average over all quenched environments,
$\mb P_\lambda^x:=\mc P\otimes P_{\omega^\lambda}^x$, is called the
\textit{annealed law}.\vspace*{1pt}
Expectations with respect to $P_{\omega^\lambda}^x$ and $\mb
P_\lambda^x$ are denoted by
$E_{\omega^\lambda}^x$ and $\mb E_\lambda^x$, respectively.\vspace*{1.5pt} We omit
the superscript when $x$ is the origin $o:=(0,\ldots, 0)$, for
example, we write $P^o_\omega$ as $P_\omega$.
%
We define the \textit{local drift} of a function $f\dvtx \mb{Z}^d\times\{e\in
\mb{Z}^d\dvtx  \llvert e\rrvert=1\}\to\mb R$ by
\[
d(f):=\sum_{e\dvtx \llvert e\rrvert=1}f(o,e)e
\]
and its spatial shift $\theta^x f$ as
\[
\theta^x f(y, e):=f(x+y,e).\vadjust{\goodbreak}
\]
When the original environment $\omega$ is deterministic and
homogeneous (i.e., $\omega=\theta^x\omega$, $\forall x$),
Sabot (\cite{Sabot}, Theorem~1) got the following perturbation expansion
for $d\ge2$:
%

\begin{quote} 
If one of $E_{\mc P}[d(\xi)]\neq o$ and $d(\omega)
\neq o$ holds, then, for $\lambda>0$
small enough, $\lim_{n\to\infty}X_n/n:=v_\lambda$
exists $\annp$-almost surely, and
\[
v_\lambda=d(\omega)+\lambda E_{\mc P}
\bigl[d(\xi)\bigr]+\lambda^2d_2+o\bigl(\lambda
^{3-\varepsilon}\bigr)\qquad \forall\varepsilon>0.
\]
The constant $d_2$ can be expressed in terms of
the Green function.
\end{quote}

%
\noindent (Sabot also obtained the expansion for $d=1$, with $d_2$ replaced by
$d_{2,\lambda}$. But in this case $v_\lambda$ can be explicitly
computed, and hence is not as interesting. See remarks in \cite{Sabot},
page~2999.)
Note that the condition for the above expansion is essentially that
$\omega^\lambda$ is \textit{ballistic} for all small $\lambda>0$, that is,
$\lim_{n\to\infty}X_n/n\neq0$ is a deterministic constant, $\annp$-a.s.

The purpose of our article is to generate Sabot's first-order expansion
to the case where the original environment is random. For RWRE in
$\mathbb Z^d, d\ge2$, two notable ballisticity conditions are
Kalikow's condition and Sznitman's \textup{(T$'$)} condition, which are introduced
in \cite{Ka81} and \cite{Sz3}, respectively. We recall that the \textup{(T$'$)}
condition is conjectured to be equivalent to the ballisticity of RWRE,
and it implies Kalikow's condition.
In this paper we are interested in two cases:
\begin{longlist}[(ii)]
\item[(i)] The original environment has zero drift (or \textit{balanced}), and $(\omega,xi)$ satisfies a Kalikow-type condition for
small $\lambda>0$:
for some $\ell\in S^{d-1}$,
{\renewcommand{\theequation}{K}
\begin{equation}\label{eqK}
\inf_{f\in\mc F}E_{\mc P} \biggl[\frac{d(\xi)\cdot\ell}{\sum
_{e\dvtx \llvert e\rrvert=1}\omega(o,e)f(e)} \biggr]
\Big/ E_{\mc P} \biggl[\frac{1}{\sum_{e\dvtx \llvert e\rrvert=1}\omega
(o,e)f(e)} \biggr]>0,
\end{equation}}\setcounter{equation}{0}%
$\mc F$ denotes the collection of nonzero functions $f\dvtx \{e\dvtx \llvert
e\rrvert=1\}\to[0,1]$.
\item[(ii)] The original environment satisfies Sznitman's ballisticity
condition \textup{(T$'$)}.
\end{longlist}
Condition (\ref{eqK}) guarantees that $\omega^\lambda$ has a speed of size
$\sim c\lambda$. Note that it is satisfied for some interesting cases,
for example, it holds for a perturbation that is ``either neutral or
pointing to the right'' (see Remark~\ref{rm1}). For the definition of
Sznitman's \textup{(T$'$)} condition, we refer to equation~(0.5) in \cite{Sz3}.

\subsection{Results}
Before the statement of our results, let us recall that one of the main
tools in the study of RWRE is the environment viewed from the point of
view of the particle process $(\bar\zeta_n)_{n\in\mb N}$, which is
defined as
\[
\bar\zeta_n=(\bar\omega_n,\bar\xi_n):=
\theta^{X_n}\zeta, \qquad n\in\mb N.
\]
Lawler \cite{Lawler} proved that for balanced environment, there
exists an ergodic invariant measure for $(\bar\zeta_n)$ which is
absolutely continuous with respect to $\mc P$.
For ballistic environment whose regeneration time has finite moment
(e.g., an environment that satisfies Sznitman's condition), it is shown
in \cite{SZ}, Theorem~3.1, that the law of~$\bar\zeta_n$ converges
weakly to an invariant measure. Recently, Berger, Cohen and Rosenthal
\cite{BCR} proved that for dimensions $d\ge4$, this measure is
ergodic and absolutely continuous with respect to the original law of
the environment.\vadjust{\goodbreak}

We denote by $\mc Q$ (for both the balanced and the ballistic cases)
the invariant measure of $(\bar\zeta_n)$ viewed from the original
RWRE, and by $\mc Q_\lambda$ the invariant measure of $(\bar\zeta
_n)$ viewed from the perturbed RWRE.

Our main results are the following.

\begin{theorem}\label{ER0}
Assume that the original environment is balanced [i.e., \mbox{$d(\omega)=o$}
almost surely] and $\mc P$ satisfies (\ref{eqK}), then
\[
\mc Q_\lambda\Rightarrow\mc Q\qquad\mbox{as }\lambda\to0,
\]
where $\Rightarrow$ denotes weak convergence.
\end{theorem}

%
\begin{theorem}\label{thm6}
Assume the $\mc P$-law of $\omega$ satisfies Sznitman's condition
\textup{(T$'$)}. Then, there exists a linear operator $\Lambda$ such that
\[
\lim_{\lambda\to0}\frac{\mc Q_\lambda f-\mc Q f}{\lambda}=\Lambda f
\]
for all (a.s.) bounded $f\dvtx \Omega\to\mb{R}$ which is
%
\begin{equation}
\label{f} \sigma\bigl(\omega_x\dvtx \llvert x \rrvert<
N_f\bigr)\mbox{-measurable for some constant }N_f
\ge1.
\end{equation}
Here, $\mc Q f$ denotes the expectation of $f$ under $\mc Q$.
Moreover, $\Lambda$ can be expressed in terms of the regeneration
times; see (\ref{e24}).
\end{theorem}

As a corollary of the above theorems, we obtain the following.

\begin{corollary}\label{thm1}
If either \textup{(i)} or \textup{(ii)} is satisfied, then
there is $\lambda_0\in[0,\kappa/2)$ such that for $\lambda\in
(0,\lambda_0)$, the limit
\[
\lim_{n\to\infty}\frac{X_n}{n}=:v_\lambda
\]
exists $\annp$-almost surely and [for the convenience of the notation,
we set $\Lambda\equiv0$ when $\mc P$ satisfies \textup{(i)}]:
\[
\lim_{\lambda\to0}\frac{v_\lambda-v_0}{\lambda} =\mc Q\bigl(d(\xi)\bigr
)+\Lambda
\bigl(d(\omega) \bigr).
\]
\end{corollary}

Recalling that for random walks in balanced random environment, Lawler
\cite{Lawler} proved that the scaling limit of $X_{\lfloor\cdot
n\rfloor}/\sqrt n$ converges to a Brownian motion with diffusion matrix
\[
D:=\bigl(E_{\mc Q}\bigl[2\omega(o,e_i)\bigr]
\delta_{ij}\bigr)_{1\le i, j\le d},
\]
the Einstein relation of a balanced random environment is an immediate
consequence of Corollary~\ref{thm1}.

\begin{proposition}[(Einstein relation)]\label{prop7}
Assume that $\mc P$-almost surely, the original environment is
balanced, and
%
\begin{equation}
\label{e56} \xi(x,e)=\omega(x,e)e\cdot\ell\qquad\forall x, e.\vadjust{\goodbreak}
\end{equation}
Then $\annp$-almost surely, $v_\lambda:=\lim_{n\to\infty
}X_n/n=\lambda E_{Q_\lambda}[d(\omega)]$, and
%
\begin{equation}
\label{ER} \lim_{\lambda\to0}\frac{v_\lambda}{\lambda}=D\ell.
\end{equation}
\end{proposition}

The zero-drift case (Theorem~\ref{ER0}) is more delicate and makes the
main part of the paper. Its proof consists of proving the following two
theorems.

\begin{theorem}\label{ER1}
Assume that the original environment is balanced. Then, for $\mc
P$-almost every $\zeta$ and any bounded measurable function $f\dvtx \Omega
\to\mathbb R$,
\[
\lim_{\lambda\to0}\frac{\lambda^2}{t} E_{\omega^\lambda} \Biggl[\sum
_{i=0}^{\lceil t/\lambda^2\rceil}f(\bar\zeta_i)
\Biggr]=\mc Q f\qquad\forall t>0.
\]
\end{theorem}

%
\begin{theorem}\label{ER2}
Assume that the original environment is balanced and $\mc P$ satisfies (\ref{eqK}).
Then for any $f$ that satisfies (\ref{f}),
\[
\Biggl\llvert\mc Q_\lambda f-\frac{\lambda^2}{t}\mb{E}_\lambda
\sum_{i=0}^{\lceil t/\lambda^2\rceil}f(\bar\zeta_i)
\Biggr\rrvert\le\frac{C\llVert f\rrVert_\infty}{\sqrt t}
\]
for all $\lambda\in(0,1/N_f)$ and $t>0$.
\end{theorem}



Our proof of Theorem \ref{ER1} is an adaption of the argument of
Lebowitz and Rost \cite{Le} (see also \cite{GMP}, Proposition 3.1) to
the discrete setting.
Namely, using a change of measure argument, we observe that the $\annp
$-law of the rescaled process $\lambda X_{\cdot/\lambda^2}$ converges
to a Brownian motion with drift.
For the proof of Theorem \ref{ER2}, we want to follow the strategy of
Gantert, Mathieu and Piatnitski \cite{GMP}. Arguments in \cite{GMP},
Proposition 5.1, show that if there is a sequence of random times
$\tau_n\sim n/\lambda^2$ (called the \textit{regeneration times}) that
divides the random path into i.i.d. parts, then good moment estimates
of the regeneration times yield the Einstein relation. [Note that the
usual definition of regeneration times, i.e., the $T(n)$'s in
Section~\ref{secball}, does not give the correct scale.] Their
definition of the regeneration times, which is a variant of that in
\cite{Shen}, crucially employs a heat kernel estimate \cite{GMP},
Lemma~5.2, for reversible diffusions.
However, due to the lack of reversibility, we do not have a heat kernel
estimate for RWRE. In this paper, we construct the regeneration times
differently, so that they divide the random path into 1-\textit{dependent}
pieces. Moreover, our regeneration times have good moment bounds, which
lead to a proof of Theorem \ref{ER2}. The key ingredients in our
construction are Kuo and Trudinger's \cite{KT} Harnack inequality for
discrete harmonic functions and the ``$\varepsilon$-coins'' trick
introduced in \cite{CZ1}.

The proof of the ballistic case (Theorem~\ref{thm6}) uses a
modification of Lebowitz and Rost's argument and the (usual)
regeneration structure for a ballistic RWRE.
The reason that the ballistic case is easier to analyze is that the
original environment already has a regeneration structure, which
provides us enough information on\vadjust{\goodbreak} the rate of the convergence to the
stationary measure. [Recall that Sznitman's \textup{(T$'$)} condition implies that
the inter-regeneration time has stretch-exponential moment.]

The structure of the paper is as follows. We will prove Theorem~\ref
{ER1} in
Section~\ref{sec1}. In Section~\ref{SecKalikow}, using Kalikow's
random walks, we obtain estimates that will be useful in deriving the
moment bounds of the regenerations.
In Section~\ref{secreg}, we present our new construction of the
regeneration times and show that they have good moment bounds.
Sections~\ref{secpro}~and~\ref{secball} are devoted
to the proofs of Theorems~\ref{ER0} and~\ref{thm6}. With these
two theorems, we obtain the derivative of the speed (w.r.t. the size of
the perturbation) in Section~\ref{secE}.

Throughout this paper, we use $c, C$ to denote finite positive
constants that depend only on the environment measure $\mc P$ (and
implicitly, on the dimension $d$ and the ellipticity constant $\kappa
$). They may differ from line to line. We also use $c_i, C_i$ to
distinguish different constants that are fixed throughout. Let $\{
e_1,\ldots, e_d\}$ be the natural basis of $\mb Z^d$.

\section{Proof of Theorem \texorpdfstring{\protect\ref{ER1}}{5}}\label{sec1}
We first consider the Radon--Nikodym derivative of the measure
$P_{\omega^\lambda}$ with respect to
$P_\omega$. For $s>0$, put
\begin{eqnarray*}
G(s, \lambda)&=&G(s,\lambda;\zeta, X_\cdot)
=\log\prod
_{j=0}^{\lceil s\rceil-1}\biggl[1+\lambda\frac{\bar\xi
_j(o,X_{j+1}-X_j)}{\bar\omega_j(o,X_{j+1}-X_j)}\biggr]
\\
&=:& \log\prod_{j=0}^{\lceil s\rceil-1}\bigl[1+\lambda a(
\bar\zeta_j, \Delta X_j)\bigr],
\end{eqnarray*}
where $\Delta X_i:=X_{i+1}-X_i$ and
\[
a(\zeta,e):=\frac{\xi(o,e)}{\omega(o,e)}.
\]
Then, for any measurable function $F$ on $C([0,s],\mathbb{R}^d)$,
\[
E_{\omega^\lambda} F(X_r\dvtx  0\le r\le s) = E_\omega
\bigl[F(X_s\dvtx  0\le r\le s)e^{G(s,\lambda)}\bigr].
\]
In particular,
%
\begin{equation}
\label{e0} E_\omega e^{G(s,\lambda)}=E_{\omega^\lambda}[1]=1
\end{equation}
for any $\lambda\in(0,1)$ and $s>0$.
Moreover, by Taylor's expansion,
%
\begin{eqnarray}
\label{e1} G(s,\lambda) &=& \sum_{j=0}^{\lceil s\rceil-1}
\log\bigl(1+\lambda a(\bar\zeta_j,\Delta X_j)\bigr)
\nonumber
\\
&=& \sum_{j=0}^{\lceil s\rceil-1} \biggl[\lambda a(\bar
\zeta_j,\Delta X_j)-\frac{\lambda^2a(\bar
\zeta_j,\Delta X_j)^2}{2} \biggr] +
\lambda^3\lceil s\rceil H
\\
&=&\lambda\sum_{j=0}^{\lceil s\rceil-1}a(\bar
\zeta_j,\Delta X_j) -\frac{\lambda^2}{2}\sum
_{j=0}^{\lceil s\rceil-1} a(\bar\zeta_j,\Delta
X_j)^2 +\lambda^3\lceil s\rceil H,\nonumber
\end{eqnarray}
where the random variable $H=H(\lambda,\zeta,X_\cdot)$ satisfies
$0\le H\le(1+\kappa^{-1})/3$.
Setting
\[
h(\zeta)=\sum_{e\dvtx \llvert e\rrvert=1}\xi(o,e)^2/
\omega(o,e),
\]
we have
\[
\Biggl( \sum_{j=0}^n \bigl[a(\bar
\zeta_j,\Delta X_j)^2-h(\bar
\zeta_j) \bigr] \Biggr)_{n\ge0}
\]
is a $P_\omega$-martingale with bounded increments. Thus, $P_\omega
$-almost surely,
\[
\lim_{n\to\infty}\frac{1}{n} \sum
_{j=0}^n \bigl[a(\bar\zeta_j,\Delta
X_j)^2-h(\bar\zeta_j) \bigr]=0.
\]
Further, recall that $\mathcal{Q}$ is the ergodic invariant measure
for $(\bar\zeta_n)_{n\ge0}$ (under $P_\omega$) and
$\mathcal{Q}\thickapprox\mathcal{P}$. Hence, by the ergodic theorem,
$\mathcal{P}\otimes P_\omega$-almost surely,
%
\begin{equation}
\label{e2} \lim_{\lambda\to0} \lambda^2\sum
_{j=0}^{\lceil t/\lambda^2\rceil-1} a(\bar\zeta_j,\Delta
X_j)^2 = \lim_{\lambda\to0}
\lambda^2\sum_{j=0}^{\lceil t/\lambda^2\rceil-1} h(\bar
\zeta_{j-1}) = tE_{\mathcal Q} h
\end{equation}
and
%
\begin{equation}
\label{e30} \lim_{\lambda\to0}\frac{\lambda^2}{t} \sum
_{i=0}^{\lceil t/\lambda^2\rceil}f(\bar\zeta_i)=E_{\mathcal Q}
f.
\end{equation}
Moreover, observing that $J_n:=\sum_{j=0}^n a(\bar\zeta_j,\Delta
X_j)$ is a $P_\omega$-martingale, by (\ref{e2}) and~\cite{DR},
Theorem~7.7.2, we get an invariance principle:

\begin{quote}
For $\mathcal{P}$-almost every $\zeta$, the process
$(\lambda J_{s/\lambda^2})_{s\ge0}$ converges weakly (under $P_\omega
$) to a Brownian motion $(N_s)_{s\ge0}$ with diffusion constant
$E_{\mathcal Q}h$.
\end{quote}

Hence, by (\ref{e1}), (\ref{e2}), (\ref{e30}) and the invariance principle,
for $\mathcal P$-almost all $\zeta$,
%
\begin{equation}
\label{e3} \frac{\lambda^2}{t} \sum_{i=0}^{\lceil t/\lambda^2\rceil}f(
\bar\zeta_i)e^{G(t/\lambda
^2,\lambda)}
\end{equation}
converges weakly to
\[
(E_{\mathcal Q}f)\exp(N_t-tE_{\mathcal Q}h/2).
\]

Next, we will prove that for $\mathcal P$-almost every $\zeta$, this
convergence is also in $L^1(P_\omega)$.
It suffices to show that the class
$(e^{G(t/\lambda^2,\lambda)})_{\lambda\in(0,1)}$ is uniformly
integrable under $P_\omega$, $\mathcal P$-a.s.
Indeed, for any $\gamma>1$, it follows from (\ref{e1}) and the
estimate on $H$ that
\begin{eqnarray*}
\gamma G\bigl(t/\lambda^2,\lambda\bigr)
&\le& G\bigl(t/\lambda^2, \gamma\lambda\bigr)
\\
&&{} + \frac{(\gamma^2-\gamma
)\lambda^2}{2}
\sum_{j=0}^{\lceil t/\lambda^2 \rceil-1}a(\bar\zeta_j,
\Delta X_j) +C\bigl(1+\gamma^3\bigr)\lambda t
\\
&<& G\bigl(t/\lambda^2, \gamma\lambda\bigr) +C\gamma^3
(t+1).
\end{eqnarray*}
Hence, for $\gamma>1$ and all $\lambda\in(0,1)$,
%
\begin{equation}
\label{e51} \quad E_\omega\exp\bigl(\gamma G\bigl(t/\lambda^2,
\lambda\bigr)\bigr) \le e^{C\gamma^3 (t+1)}E_\omega\exp\bigl(G\bigl(t/
\lambda^2, \gamma\lambda\bigr)\bigr) \stackrel{\fontsize{8.36pt}{8.36pt}{\mathrm{by}~(\ref{e0})}}
{=}e^{C\gamma^3 (t+1)},
\end{equation}
which implies the uniform integrability of $(e^{G(t/\lambda^2,\lambda
)})_{\lambda\in(0,1)}$.
So the $E_\omega$-expecta\-tion of (\ref{e3}) also converges to the
expectation of its weak limit (for $\mathcal P$-almost every~$\zeta$) and
\[
\lim_{\lambda\to0} E_{\omega^\lambda} \Biggl[\frac{\lambda^2}{t} \sum
_{i=0}^{\lceil t/\lambda^2\rceil}f(\bar\zeta_i)
\Biggr] =(E_{\mathcal Q}f)E \bigl[\exp(N_t-t E_{\mathcal Q}h/2)
\bigr].
\]
The theorem follows by noting that $tE_{\mathcal Q}h=EN_t^2$ and that
%
\[
E \bigl[\exp\bigl(N_t-EN_t^2/2\bigr)
\bigr] =1.
\]
%

\section{Kalikow's auxiliary random walks}\label{SecKalikow}
In this section, we will recall Kalikow's auxiliary random walks and
use it to obtain some estimates that will be useful later.

For any connected strict subset $U$ of $\mb Z^d$, let
\begin{eqnarray*}
\partial U &=& \bigl\{x\in\mb Z^d\setminus U\dvtx  \exists y\in U, \llvert
y-x\rrvert=1\bigr\},
\\
T_U&=&\inf\{n\ge0\dvtx  X_n\in\partial U\}.
\end{eqnarray*}
Define on $U\cup\partial U$ a Markov chain with transition probability
%
\begin{equation}
\label{e54} \hat P_U(x,x+e)= \cases{ \displaystyle
\frac{E_{\mc P}E_\omega[\sum_{n=0}^{T_U}1_{X_n=x}\omega(x,e)]}{
E_{\mc P}E_\omega[\sum_{n=0}^{T_U}1_{X_n=x}]}, &\quad$x\in U$, $\llvert
e\rrvert=1$,
\cr
1, &\quad$x\in
\partial U$, $e=o$,}
\end{equation}
and set
\[
\hat d_{U}(x):=\sum_{e}e\hat
P_U(x,x+e).
\]
We say that the \textit{Kalikow's condition relative to $\ell\in S^{d-1}$} holds if there exists $\delta>0$ such that
%
\begin{equation}
\label{Kalikow} \inf_{U,x\in U}\hat d_{U}(x)\cdot\ell\ge
\delta.
\end{equation}
The interest of this Markov chain lies in the fact that $\hat P_U$ and
$\mb P$ have the same exit distribution from $U$ (\cite{Ka81},
Proposition 1):
\[ 
\mbox{if } \hat P_U(T_U<\infty)=1,\qquad\mbox{then }
\hat P_U(X_{T_U}\in\cdot)=\mb P(X_{T_U}\in\cdot).
\] 

%
\begin{theorem}[(\cite{SZ}, Theorem 2.3)]
If (\ref{Kalikow}) holds, then there exists a deterministic $v\in\mb
R^d$ such that
\[
\lim_{n\to\infty}X_n/n=v,\qquad\mb P\mbox{-a.s.}
\]
\end{theorem}

It is also shown in \cite{Ka81}, (11), that
(\ref{Kalikow}) has the following sufficient condition:
%
\begin{equation}
\label{Kalikow3} \inf_{f\in\mc F} E_{\mc P} \biggl[
\frac{d(\omega)\cdot\ell}{\sum_e\omega(o,e)f(e)} \biggr] \Big/E_{\mc P}
\biggl[\frac{1}{\sum_e\omega(o,e)f(e)} \biggr]\ge
\delta,
\end{equation}
where $\mc F$ is the same as in (\ref{eqK}).

%
\begin{proposition}\label{prop9}
Assume \textup{(i)}.
Then for some $\lambda_0>0$ and all $\lambda\in[0,\lambda_0)$,
there is a deterministic constant $v_\lambda\in\mathbb{R}^d$ such that
\[
\lim_{t\to\infty} \frac{X_t}{t}=v_\lambda, \qquad
\mb{P}_\lambda\mbox{-almost surely}.
\]
\end{proposition}

\begin{pf}
By (\ref{eqK}), there exist $\lambda_0>0$ such that for all $\lambda\in
(0,\lambda_0)$ and
\[
\inf_{f\in\mc F} E_{\mc P} \biggl[\frac{d(\xi)\cdot\ell}{\sum
_{e\dvtx \llvert e\rrvert=1}\omega^\lambda(o,e)f(e)}
\biggr] \Big/ E_{\mc P} \biggl[\frac{1}{\sum_{e\dvtx \llvert e\rrvert
=1}\omega
^\lambda(o,e)f(e)} \biggr] >0.
\]
Noting that $d(\omega)=0$, there is $\rho>0$ such that the law of
$\omega^\lambda$ satisfies (\ref{Kalikow3}), with $\delta$ replaced
by $\lambda\rho$. This implies
%
\begin{equation}
\label{e55} \inf_{U,x\in U}\hat d_{U}(x)\cdot\ell\ge
\lambda\rho.
\end{equation}
The proposition follows.
\end{pf}

%
\begin{remark}\label{rm1}
Although (\ref{eqK}) looks complicated, it includes some simple cases:
\begin{longlist}[(a)]
\item[(a)] (\ref{eqK}) holds when
\[
E_{\mc P}\bigl[\bigl(d(\xi)\cdot\ell\bigr)_+\bigr]>\frac{1}{\varepsilon
}E_{\mc P}
\bigl[\bigl(d(\xi)\cdot\ell\bigr)_-\bigr].
\]
For instance, (\ref{eqK}) is satisfied when the perturbation is ``either
neutral or pointing to the right.'' See \cite{SZ}, Proposition 2.4.

\item[(b)] When $\omega$ and $\xi$ are independent, (\ref{eqK}) is
equivalent to $E_{\mc P}[d(\xi)]\neq o$.
\end{longlist}
\end{remark}
%
\subsection{Auxiliary estimates}
In this subsection, we consider perturbed RWRE that satisfies (i).
Making use of Kalikow's random walks, we obtain some auxiliary
estimates that will be useful in getting the regeneration moment bounds
in Section~\ref{secreg}.

From now on, we assume that (i) holds with
\[
\ell=e_1.
\]
(The same arguments work also for general $\ell\in S^{d-1}$, but with
cumbersome notations.)
Recall that (i) implies (\ref{e55}):
\[
\inf_{U,x\in U}\hat d_{U}(x)\cdot e_1\ge
\lambda\rho.
\]

Let
\[
\lambda_1:=0.5/\bigl\lceil(2\lambda)^{-1}\bigr\rceil
\]
so that $0.5/\lambda_1$ is an integer and
\[
\frac{1}{2\lambda}\le\frac{1}{2\lambda_1}<\frac{1}{2\lambda}+1.
\]
%

For any $k\in\frac{1}{2}\mathbb{Z}, x\in\mathbb{Z}^d$, set
%
\begin{eqnarray}
\label{e57} \mathcal{H}_k^x &=&\mathcal{H}_n^x(
\lambda,\ell):= \bigl\{y\in\mathbb{Z}^d\dvtx  (y-x)\cdot e_1=k/
\lambda_1\bigr\},
\\
\label{e48} T_k &:=&\inf\bigl\{t\ge0\dvtx  (X_t-X_0)
\cdot e_1=k/\lambda_1\bigr\}.
\end{eqnarray}
For $n\in\mathbb N$, we call $\mc H_n^x$ \textit{the $n$th level}
(with respect to $x$).
Since the random walk is transient in the $e_1$ direction, $T_k$'s are
finite $\mb P_\lambda$-almost surely.

\begin{proposition}\label{prop6}
Let $(X'_n)$ be a simple random walk on $\mb{Z}$ with
\[
\frac{P(X'_{i+1}=x+1\mid X_i'=x)}{P(X'_{i+1}=x-1\mid X_i'=x)}=q\neq
1\qquad\forall x\in\mb{Z}.
\]
Then for any $i,j\in\mb{Z}^+$ and $-j\le0\le i$,
\[
P\bigl(X' \mbox{ visits $-j$ before visiting }i\mid
X'_0=0\bigr)=\frac{q^i-1}{q^{i+j}-1}.
\]
In particular, when $q<1$,
\[
P\bigl(X' \mbox{ never visits }{-}j\mid X'_0=0
\bigr)=1-q^{-j}.
\]
\end{proposition}

The proof is omitted.

\begin{proposition}\label{prop5}
Assume \textup{(i)}. There exists $\lambda_0\in(0,1)$ such that for $\lambda
\in(0,\lambda_0)$ and any $n,m\in\mb{N}/2$:
\begin{longlist}[(a)]
\item[(a)] $\frac{e^{5m/\kappa}-1}{e^{5(m+n)/\kappa}-1}\le\mb
{P}_\lambda(T_{-n}<T_m)\le\frac{e^{m\rho/2}-1}{e^{(n+m)\rho/2}-1}$;\vspace*{1pt}

\item[(b)] $\frac{e^{5m/\kappa}-1}{e^{5(m+n)/\kappa}-1}\le P_{\omega
^\lambda}(T_{-n}<T_{m})\le\frac{1-e^{-5m/\kappa
}}{1-e^{-5(m+n)/\kappa}}$, $\mc P$-almost surely.\vadjust{\goodbreak}
\end{longlist}
\end{proposition}

\begin{pf}
(a) Let $U= \{z\in\mb{Z}^d\dvtx  n/\lambda_1\le z\cdot e_1\le
m/\lambda_1\}$. With abuse of notation, we let $\hat{X}$ be the
Markov chain defined at (\ref{e54}), with $\omega$ replaced by
$\omega^\lambda$ (because we are interested in the perturbed
environment). Since [by (\ref{e55}), (\ref{e54}) and $\lambda<\kappa/2$]
\[
1+\rho\lambda\le\frac{\hat P_U(x,x+e_1)}{\hat P_U(x,x-e_1)}\le1+\frac
{4\lambda}{\kappa},
\]
we can couple two Markov chains $X', X''$ on $\mb{Z}$ to $\hat X$ such
that for all $i\in\mb N$, $x\in\mb Z$,
\begin{eqnarray*}
\frac{P(X'_{i+1} = x+1\mid X'_i=x)}{P(X'_{i+1}=x-1\mid X'_i=x)}&=&1+\rho
\lambda,
\\
\frac{P(X''_{i+1}=x+1\mid X''_i=x)}{P(X''_{i+1}=x-1\mid
X''_i=x)}&=&1+\frac{4\lambda}{\kappa},
\end{eqnarray*}
and
\[
X'_i\le\hat X_i\cdot e_1\le
X''_i\qquad\forall i\in\mb N.
\]
Hence, by Proposition~\ref{prop6}, we obtain
\[
\frac{(1+4\lambda/\kappa)^{m/\lambda_1}-1}{(1+4\lambda/\kappa
)^{(n+m)/\lambda_1}-1}\le\hat{P}_U(T_{-n}<T_m)
\le\frac{(1+\rho
\lambda)^{m/\lambda_1}-1}{(1+\rho\lambda)^{(n+m)/\lambda_1}-1}.
\]
Taking $\lambda$ small enough, inequality (a) is proved.

(b)~Observe that for $\lambda\in(0,\kappa/2)$, $\mc P$-almost surely,
\[
1-\frac{4\lambda}{\kappa}\le\frac{\omega^\lambda(x,e_1)}{\omega
^\lambda(x,-e_1)}\le1+\frac{4\lambda}{\kappa}.
\]
Inequality (b) then follows from the same argument as in the proof of (a).
\end{pf}

%
\begin{theorem}\label{44444}
Assume that $\omega$ is balanced. Let
\[
\tilde T_n:=T_n\wedge T_{-n}.
\]
There exists a constant $s>0$ such that for any uniformly elliptic
balanced environment $\omega$ and all $\lambda\in(0,\kappa/2), n\in
\mb N$,
\[
E_{\omega^\lambda}\bigl[e^{s\lambda^2\tilde T_n/n}\bigr]<C.
\]
\end{theorem}

The proof, which uses coupling, is given in the \hyperref[appe]{Appendix}.

%
\begin{proposition}\label{prop8}
Assume that $\omega$ is balanced. There exists a constant $C_0$ such
that for $\mc P$-almost all $(\omega,\xi)$ and $\lambda\in(0,\kappa/2)$,
\[
P_{\omega^\lambda}\bigl(\llvert X_{T_{0.5}}\rrvert<C_0/
\lambda\bigr)>1/C_0.
\]
\end{proposition}

\begin{pf}
Let $Y_n:=X_n-\lambda\sum_{i=0}^{n-1}d(\theta^{X_i} \xi)$. Then
$(Y_n)$ is a $P_{\omega^\lambda}$-martingale.
Recall the definition of $\tilde T_n$ in Theorem~\ref{44444}. For any
$K>0$ and $\tilde K:=K/(4d)$,
\begin{eqnarray*}
&& P_{\omega^\lambda}\bigl(\llvert X_{T_{0.5}}\rrvert<K/\lambda\bigr)
\\
&&\qquad \ge1-P_{\omega^\lambda}\Bigl(\max_{t\le\tilde K/\lambda^2}\llvert
X_t\rrvert\ge K/\lambda\Bigr)-P_{\omega^\lambda}\bigl(\tilde
T_{0.5}>\tilde K/\lambda^2\bigr) -P_{\omega^\lambda}(T_{0.5}>T_{-0.5}).
\end{eqnarray*}
By Proposition~\ref{prop5}(b) and Theorem~\ref{44444}, it suffices to
show that
\[
P_{\omega^\lambda}\Bigl(\max_{t\le\tilde K/\lambda^2}\llvert
X_t\rrvert\ge K/\lambda\Bigr)
\]
can be sufficiently small if $K$
is large. Indeed,
\begin{eqnarray*}
P_{\omega^\lambda}\Bigl(\max_{t\le\tilde K/\lambda^2}\llvert X_t
\rrvert\ge K/\lambda\Bigr) &\le& P_{\omega^\lambda}\Bigl(\max_{t\le
\tilde K/\lambda^2}
\llvert Y_t\rrvert\ge CK/\lambda\Bigr)
\\
&\le& C e^{-(c(K/\lambda)^2)/(K/\lambda^2)}=Ce^{-cK},
\end{eqnarray*}
where we used Azuma--Hoeffding inequality in the last inequality.
\end{pf}

%
\begin{lemma}\label{prop2}
Assume \textup{(i)}; then $\mb P_\lambda(\llvert X_{T_n}-\frac{n}{\lambda
_1}e_1\rrvert\ge\frac{n}{\lambda_1})\le Ce^{-cn}$.
\end{lemma}

\begin{pf}
Observe that
\[
\mb P_\lambda\biggl(\biggl\llvert X_{T_n}-
\frac{n}{\lambda_1}e_1\biggr\rrvert\ge\frac{n}{\lambda_1}\biggr) = \hat
P_{U_n}\biggl(\biggl\llvert\hat X_{T_n}-
\frac{n}{\lambda_1}e_1\biggr\rrvert\ge\frac{n}{\lambda_1}\biggr),
\]
where $U_n= \{x\dvtx x\cdot e_1\le n/\lambda_1\}$.
For $j\ge0$, let
\[
\hat Y_j:=\hat X_j-\sum
_{i=0}^{j-1}\hat d_{U_n}(X_i).
\]
Then, for $k_n:=\frac{2n}{\rho\lambda^2}$,
%
\begin{eqnarray}
\label{e45}
&& \hat P_{U_n} \biggl(\biggl\llvert\hat
X_{T_n}-\frac{n}{\lambda
_1}e_1\biggr\rrvert\ge
\frac{5n}{\rho\lambda_1} \biggr)
\nonumber
\\
&&\qquad = \hat P_{U_n}(T_n\ge k_n)+ \hat
P_{U_n} \biggl(\max_{0\le i\le k_n}\llvert\hat
X_i\rrvert\ge\frac{5n}{\rho\lambda_1} \biggr)
\nonumber
\\
&&\qquad \le\hat P_{U_n} \Biggl(\hat Y_{k_n}\cdot e_1+
\sum_{i=0}^{k_n-1}\hat d_{U_n}(X_i)
\cdot e_1\le\frac{n}{\lambda_1} \Biggr)
\\
&&\quad\qquad{}+\hat P_{U_n} \Biggl(
\max_{0\le i\le k_n}\llvert\hat Y_i\rrvert\ge
\frac{5n}{\rho\lambda_1}-\Biggl\llvert\sum_{i=0}^{k_n-1}
\hat d_{U_n}(X_i)\Biggr\rrvert\Biggr)
\nonumber
\\
&&\qquad \le\hat P_{U_n}\biggl(\hat Y_{k_n}\cdot e_1\le-
\frac{n}{\lambda_1}\biggr) +\hat P_{U_n}\biggl(\max_{0\le i\le k_n}
\llvert\hat Y_i\rrvert\ge\frac{n}{\rho\lambda_1}\biggr),\nonumber
\end{eqnarray}
where in the last inequality we used the fact that
\[
\rho\lambda\le\hat d_{U_n}\le2\lambda.
\]
The lemma follows by observing that $(\hat Y_j)_{j\ge0}$ is a
martingale with bounded increments and by applying the Azuma--Hoeffding
inequality to (\ref{e45}).
\end{pf}

\section{Regenerations}\label{secreg}
In this section, we will construct a 1-dependent regeneration structure
for perturbed RWRE that satisfies (i). Recall that we assume (without
loss of generality) that $\ell=e_1$.
\subsection{Harnack inequality and its application}\label{secharnack}
Let $a$ be a nonnegative function on $\mathbb{Z}^d\times\mathbb{Z}^d$
such that for any $x$,
$a(x,y)> 0$
only if $x$ and $y$ are neighbors, that is, $\llvert x-y\rrvert=1$,
denoted $x\sim y$.
We also assume that
\[
\sum_y a(x,y)=1\qquad\forall x\in
\mathbb{Z}^d.
\]
Define the linear operator $L_a$ acting on the set of functions
on $\mathbb{Z}^d$ by
\[
L_a f(x)=\sum_y a(x,y)
\bigl(f(y)-f(x)\bigr).
\]
Set
\[
b(x)=\sum_y a(x,y) (y-x)\quad\mbox{and}\quad
b_0=\sup\llvert b\rrvert.
\]
We assume that $L_a$ is uniformly elliptic with constant $\kappa\in
(0,\frac{1}{2d}]$. That is,
\[
a(x,y)\ge\kappa\qquad\mbox{for any }x, y\mbox{ such that }x\sim y.
\]
For $r>0, x\in\mathbb{R}^d$, let $B_r(x)=\{z\in\mathbb{Z}^d\dvtx  \llvert
z-x\rrvert<r\}$. We also write
$B_r(o)$ as $B_r$.

The following Harnack estimate is due to Kuo and Trudinger \cite{KT},
Theorem 3.1. See also the Appendix of \cite{Guo} for a detailed proof.

\begin{theorem}[(Harnack inequality)]\label{harnack}
Let $u$ be a
nonnegative function on $B_R$, $R>1$. If
\[
L_a u=0
\]
in $B_R$,
then for any $\sigma\in(0,1)$ with $R(1-\sigma)>1$, we have
\[
\max_{B_{\sigma R}}u\le C\min_{B_{\sigma R}}u,
\]
where $C$ is a positive constant depending on $d, \kappa, \sigma$ and
$b_0 R$.
\end{theorem}

With the Harnack inequality, we have
the following.

\begin{lemma}\label{l4}
Assume \textup{(i)}. There exists a constant $c_1\in(0,1]$ such that for
$\lambda\in(0,1)$, $x\in\mathbb{Z}^d$
and $\mc P$-almost every $(\omega,\xi)$,
%
\begin{equation}
\label{e5} P_{\omega^\lambda}^{x}(X_{T_1}=\cdot) \ge
c_1 P_{\omega^\lambda}^{x+0.5e_1/\lambda_1}(X_{T_{0.5}}=\cdot\mid
T_{0.5}<T_{-0.5}).
\end{equation}
\end{lemma}

\begin{pf}
%
For any $x\in\mathbb Z^d$ and $k\in\frac{1}{2}\mathbb Z$, recall the
definition of $\mc H^x_k$ in (\ref{e57}).
Fix $w\in\mathcal{H}_1^x$. Then the function
\[
f(z):= P_{\omega^\lambda}^z\bigl(X_\cdot\mbox{ visits }
\mathcal{H}_1^x\mbox{ for the first time at }w\bigr)
\]
satisfies
\[
L_{\omega^\lambda}f(z)=0
\]
for all $z\in\{y\dvtx  (y-x)\cdot e_1<1/\lambda_1\}$.
By Theorem \ref{harnack} (in this case $a=\omega^\lambda,
R=0.5/\lambda_1$ and $b_0\le\lambda$), there exists a constant $C_2$
such that
for any $y, z\in\mathcal{H}_{0.5}^x$ with $\llvert z-y\rrvert\le
0.5/\lambda_1$,
%
\begin{equation}
\label{e49} f(z)\ge C_2 f(y).
\end{equation}
Hence, for any $z\in\mathcal{H}_{0.5}^x$ such that $\llvert
z-(x+0.5e_1/\lambda_1)\rrvert<C_0/\lambda_1$ (recall that $C_0$
is the constant in Proposition~\ref{prop8}), we have
%
\begin{equation}
\label{e27} f(z)\ge C_2^{2C_0} f(x+0.5e_1/
\lambda_1).
\end{equation}
Therefore,
\begin{eqnarray*}
P_{\omega^\lambda}^x(X_{T_1}=w) &\ge& \sum
_{\llvert y-x\rrvert<C_0/\lambda} P_{\omega^\lambda
}^x(X_{T_{0.5}}=y)P_{\omega^\lambda
}^y(X_{T_{0.5}}=w)
\\
&\stackrel{\fontsize{8.36pt}{8.36pt}{(\ref{e27})}} {\ge}& C P_{\omega^\lambda}^x\bigl(
\llvert X_{T_{0.5}}-x\rrvert<C_0/\lambda_1
\bigr)P_{\omega^\lambda}^{x+0.5e_1/\lambda_1} (X_{T_{0.5}}=w)
\\
&\ge& c_1 P_{\omega^\lambda}^{x+0.5e_1/\lambda_1} (X_{T_{0.5}}=w\mid
T_{0.5}<T_{-0.5}),
\end{eqnarray*}
where in the last inequality we used Proposition~\ref{prop8}
and [Proposition~\ref{prop5}(b)]
\[
P_{\omega^\lambda}^{x+0.5e_1/\lambda_1}(T_{0.5}<T_{-0.5})>C.
\]\upqed
\end{pf}

\subsection{Construction of the regeneration times}
In this subsection, we will construct regeneration times that allow the
path to backtrack at most distance $1/\lambda$ in direction $e_1$
after each regeneration. The main difficulty is to decouple the parts
before and after a regeneration in such a way that they are ``almost
independent.'' Our main observation is that (by Lemma~\ref{l4}) the
hitting probability $P_{\omega^\lambda}^{x}(X_{T_1}=\cdot)$ to the
next level dominates [in the sense of (\ref{e5})] a ``good''
probability measure
%
\begin{equation}
\label{e17} \mu_{\omega^\lambda,1}^x(\cdot):=P_{\omega^\lambda
}^{x+0.5e_1/\lambda_1}
(X_{T_{0.5}}=\cdot\mid T_{0.5}<T_{-0.5}),
\end{equation}
which is independent of environment to the left of level $\mc H^x_0$.
Hence, the hitting probability can be decomposed as
\[
P_{\omega^\lambda}^x(X_{T_1}=\cdot) =\beta
\mu_{\omega^\lambda,1}^x(\cdot)+(1-\beta)\mu_{\omega
^\lambda,0}^x(
\cdot),
\]
where (recall that $c_1$ is the constant in Lemma~\ref{l4})
\[
\beta:=c_1/2 \quad\mbox{and}\quad\mu_{\omega^\lambda,0}^x(
\cdot) 
:= \bigl[P_{\omega^\lambda}^x(X_{T_1}=
\cdot)-\beta\mu_{\omega
^\lambda,1}^x(\cdot) \bigr]/(1-\beta).
\]
Note that by (\ref{e5}), both $\mu_{\omega^\lambda,1}^x$ and $\mu
_{\omega^\lambda,0}^x$ are probability measures on $\mathcal{H}^x_1$.
This suggests us to use a coin-tossing trick to decouple the paths and
define the regenerations, which we explain as follows.

\begin{figure}

\includegraphics{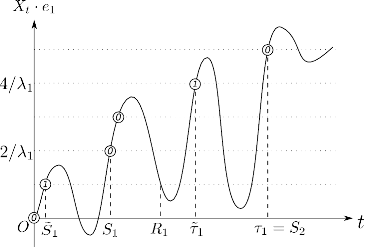}

\caption{In this picture, $K=2, X_{\tau_1}=5/\lambda_1,
M_1=4/\lambda_1$.}\label{fig1}
\end{figure}

For any $\mathcal{O}\in\sigma(X_1,X_2,\ldots, X_{T_1}), x\in
\mathbb{Z}^d$ and $i\in\{0,1\}$, put
%
\begin{equation}
\label{e6} \nu_{\omega^\lambda,i}^x(\mathcal{O})
:= \sum_y \bigl[i
\mu_{\omega^\lambda,1}^x(y)+(1-i)\mu_{\omega^\lambda
,0}^x(y)
\bigr] P_{\omega^\lambda}^x(\mathcal{O}\mid X_{T_1}=y).
\end{equation}
Let $(\varepsilon_i)_{i=1}^\infty\in\{0,1\}^\mathbb{N}$ be i.i.d.
Bernoulli random variables with law $Q_\beta$:
\[
Q_\beta(\varepsilon_i=1)=\beta\quad\mbox{and}\quad
Q_\beta(\varepsilon_i=0)=1-\beta.
\]
Intuitively, whenever the walker visits a new level $\mathcal{H}_i,
i\ge0$,
we make him flip a coin $\varepsilon_i$.
If $\varepsilon_i=0$ (or $1$), he then walks
following the law $\nu_{\omega^\lambda,0}$ (or $\nu_{\omega
^\lambda,1}$) until he reaches the $(i+1)$th level. The regeneration
time $\tau_1$ is defined to be the first time of visiting a new level
$\mathcal{H}_k$ such that the outcome $\varepsilon_{k-1}$ of the previous
coin-tossing is ``$1$'' and the path will never backtrack to level
$\mathcal{H}_{k-1}$ in the
future. See Figure~\ref{fig1}.

%

%
%
\begin{figure}[b]

\includegraphics{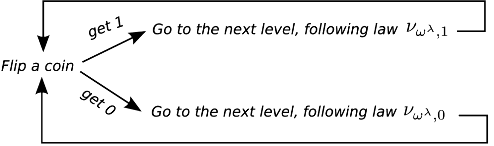}

\caption{The law $\bar P_{\omega^\lambda,\varepsilon}$ for the
walks.}\label{fig0}
\end{figure}


We now give the formal definition of the regeneration times.

We sample the sequence $\varepsilon:=(\varepsilon_i)_{i=1}^\infty$
according to the product measure $Q_\beta$ and fix it. Then we define
a new law $P_{\omega^\lambda,\varepsilon}$
on the paths, by the following steps (see Figure~\ref{fig0}):
\begin{itemize}
\item \textit{Step} 1. For $x\in\mathbb{Z}^d$, set
\[
P_{\omega^\lambda,\varepsilon}^x(X_0=x)=1.
\]
\item \textit{Step} 2.
Suppose the $P_{\omega^\lambda,\varepsilon}^x$-law for paths of
length${}\le n$ is defined. For any path $(x_i)_{i=0}^{n+1}$ with
$x_0=x$, define
\begin{eqnarray*}
&& P_{\omega^\lambda,\varepsilon}^{x} (X_{n+1}=x_{n+1},\ldots,
X_{0}=x_0)
\\
&&\qquad := P_{\omega,\varepsilon}^{x}(X_I=x_I,\ldots,
X_0=x_0) \nu_{\omega^\lambda,\varepsilon_J}^{x_I}(X_{n+1-I}=x_{n+1},
\ldots, X_1=x_{I+1}),
\end{eqnarray*}
where
\[
J=\max\bigl\{j\ge0\dvtx  \mathcal{H}_{j}^{x_0}\cap
\{x_i, 0\le i\le n\}\neq\varnothing\bigr\}
\]
is the highest level visited by $(x_i)_{i=0}^{n}$ and
\[
I=\min\bigl\{0\le i\le n\dvtx  x_i\in\mathcal{H}_J^{x_0}
\bigr\}
\]
is the hitting time to the $J$th level.
\item \textit{Step} 3.
By induction, the law $P_{\omega^\lambda,\varepsilon}^x$
is well defined for paths of all lengths.
\end{itemize}

Note that a path sampled by $P_{\omega^\lambda,\varepsilon}^x$ is not a
Markov chain, but
the law of $X_\cdot$ under
\[
\bar{P}_{\omega^\lambda}^x 
:=Q_\beta\otimes
P_{\omega^\lambda,\varepsilon}^x
\]
coincides with $P_{\omega^\lambda}^x$. That is,
%
\begin{equation}
\label{e47} \bar{P}_{\omega^\lambda}^x(X_\cdot\in\cdot) =
P_{\omega^\lambda}^x(X_\cdot\in\cdot).
\end{equation}
We denote by
$\bar{\mathbb P}_\lambda
:=\mc P\otimes\bar{P}_{\omega^\lambda}$ the law of the triple
$(\omega,\varepsilon, X_\cdot)$. Expectations with respect to $\bar
{P}_{\omega^\lambda}^x$ and $\bar{\mathbb P}_\lambda$ are denoted by
$\bar{E}_{\omega^\lambda}^x$ and $\bar{\mathbb E}_\lambda$, respectively.

\begin{figure}

\includegraphics{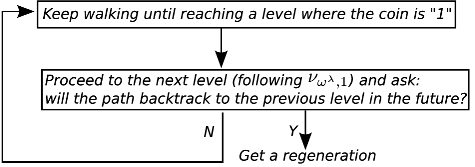}

\caption{The definition of a regeneration time.}\label{fig2}
\end{figure}

Next, for a path $(X_n)_{n\ge0}$ sampled according to $P_{\omega
^\lambda,\varepsilon}^o$, we will define the regeneration times. See
Figure~\ref{fig2} for an illustration.
%

To be specific, put $S_0=0, M_0=0$,
and define inductively
\begin{eqnarray*}
S_{k+1}&=&\inf\{T_{n+1}\dvtx  n/\lambda_1\ge
M_k\mbox{ and }\varepsilon_n=1\},
\\
R_{k+1}&=&S_{k+1}+T_{-1}\circ\theta_{S_{k+1}},
\\
M_{k+1}&=&X_{S_{k+1}}\cdot e_1+N\circ
\theta_{S_{k+1}}/\lambda_1, \qquad k\ge0.
\end{eqnarray*}
Here, $\theta_n$ denotes the time shift of the path, that is, $\theta
_n X_\cdot=(X_{n+i})_{i=0}^\infty$,
and
%
\begin{equation}
\label{e46} N:=\inf\bigl\{n\dvtx  n/\lambda_1>(X_i-X_0)
\cdot e_1\mbox{ for all }i\le T_{-1}\bigr\}.
\end{equation}
Set
\begin{eqnarray*}
K&:=&\inf\{k\ge1\dvtx  S_k<\infty, R_k=\infty\},
\\
\tau_1&:=&S_K\quad\mbox{and}\quad\tau_{k+1}=
\tau_k+\tau_1\circ\theta_{\tau_k}.
\end{eqnarray*}
We call $(\tau_k)_{k\ge1}$ \textit{regeneration times}.

\subsection{The renewal property of the regenerations}

The regeneration times possess good renewal properties:
\begin{longlist}[2.]
%
\item[1.]
Set $\tau_0=0$. For $k\ge0$, define
\begin{eqnarray*}
\tilde{S}_{k+1}&:=&\inf\{T_n\dvtx  n/\lambda\ge M_k
\mbox{ and }\varepsilon_n=1\},
\\
\tilde{\tau}_1&:=&\tilde{S}_K\quad\mbox{and set}\quad
\tilde{\tau}_{k+1}:=\tau_k+\tilde{\tau}_1\circ
\theta_{\tau_k}.
\end{eqnarray*}
Namely, $\tilde\tau_k$ is the hitting time to the previous level of
$X_{\tau_k}$.
Conditioning on \mbox{$X_{\tilde{\tau}_k}=x$}, the law of $X_{\tau_k}$ is
$\mu_{\omega^\lambda,1}^x$, which is independent (under the
environment measure $\mc P$) of
$\sigma(\zeta_y\dvtx y\cdot e_1\le x\cdot e_1)$. Moreover, after time
$\tau_k$,
the path will never visit $\{y\dvtx y\cdot e_1\le x\cdot e_1\}$. Therefore,
$\tau_{k+1}-\tau_k$ is independent of what happened before $\tau
_{k-1}$ and the inter-regeneration times form a $1$-dependent sequence.

\item[2.]
Since $(X_{\tilde\tau_{k+1}}-X_{\tau_k})_{k\ge1}$ are i.i.d. and
$(X_{\tau_{k+1}}-X_{\tilde\tau_{k+1}})\cdot e_1=1/\lambda_1$, the
inter-regeneration distances
$ ((X_{\tau_{k+1}}-X_{\tau_k})\cdot e_1 )_{k\ge1}$ are i.i.d.

\item[3.]
From the construction, we see that a regeneration occurs after roughly
a geometric number of levels. Thus, we expect $(X_{\tau_{k+1}}-X_{\tau
_k})\cdot e_1\sim c/\lambda$ and (by Theorem~\ref{44444}) $\tau
_{k+1}-\tau_k\sim c/\lambda^2$.
\end{longlist}

The above properties will be verified in Lemma~\ref{l5},
Proposition~\ref{prop1} and Corollary~\ref{cor1}.

We introduce the $\sigma$-field
\[
\mathcal{G}_k:= \sigma\bigl( \tilde{\tau}_k,
(X_i)_{i\le\tilde{\tau}_k},(\zeta_y)_{y\cdot
e_1\le X_{\tilde{\tau}_k}\cdot e_1} \bigr)
\]
and set
%
\begin{equation}
\label{e18} p_\lambda:=E_{\mc P} \biggl[\sum
_y\mu_{\omega^\lambda,1}(y)P_{\omega
^\lambda}^y(T_{-1}=
\infty) \biggr].
\end{equation}
%

\begin{lemma}\label{l5}
For any appropriate measurable sets $B_1, B_2$
and any event
\[
B:=\bigl\{(X_i)_{i\ge0}\in B_1, (
\zeta_y)_{y\cdot e_1>-1/\lambda_1}\in B_2\bigr\},
\]
we have for $k\ge1$,
\[
\bar{\mathbb P}_\lambda(B\circ\bar{\theta}_{\tau_k}\mid
\mathcal{G}_k) = E_{\mc P} \biggl[ \sum
_y \mu_{\omega^\lambda,1}(y) \bar{P}_{\omega^\lambda}^y
\bigl(B\cap\{T_{-1}=\infty\}\bigr) \biggr] \Big/p_\lambda,
\]
where $\bar{\theta}_n$ is the time-shift defined by
\[
B\circ\bar{\theta}_n = \bigl\{(X_i)_{i\ge n}\in
B_1, (\zeta_y)_{(y-X_n)\cdot e_1>-1/\lambda_1}\in B_2\bigr
\}.
\]
\end{lemma}

\begin{pf}
First, we consider the case $k=1$.
Let $\vartheta^n$ denote the shift of the \mbox{$\varepsilon$-}coins, that is,
$\vartheta^n \varepsilon_\cdot=(\varepsilon_i)_{i\ge n}$.
For any $A\in\mathcal{G}_1$,
\begin{eqnarray*}
&& \bar{\mathbb P}_\lambda(B\circ\bar{\theta}_{\tau_1}\cap A)
\\
&&\qquad = E_{\mc P\otimes Q_\beta} \biggl[ \sum_{k\ge1,x}P_{\omega^\lambda
,\varepsilon}
\bigl(A\cap\{\tilde{S}_k<\infty,R_k=\infty,
X_{\tilde{S}_k}=x\}\cap B\circ\bar{\theta}_{S_k}\bigr) \biggr]
\\
&&\qquad = E_{\mc P\otimes Q_\beta} \biggl[ \sum_{k\ge1,x,y}P_{\omega^\lambda
,\varepsilon}
\bigl(A\cap\{\tilde{S}_k<\infty,X_{\tilde{S}_k}=x\}\bigr)
\nu_{\omega^\lambda,1}^x(X_{T_1}=x+y)
\\
&&\hspace*{194pt}{} \times P_{\omega^\lambda,\vartheta
^{k+1}\varepsilon}^{x+y}\bigl(B
\cap\{ T_{-1}=\infty\}\bigr) \biggr].
\end{eqnarray*}

Note that in the last equality,
\[
P_{\omega^\lambda,\varepsilon} \bigl(A\cap\{\tilde{S}_k<\infty,X_{\tilde{S}_k}=x
\}\bigr)
\]
is
$\sigma((\varepsilon_i)_{i\le k},(\zeta_z)_{(z-x)\cdot e_1\le
0} )$-measurable,
whereas
\[
\nu_{\omega^\lambda,1}^x(X_{T_1}=x+y)P_{\omega^\lambda,\vartheta
^{k+1}\varepsilon}^{x+y}
\bigl(B\cap\{T_{-1}=\infty\}\bigr)
\]
is
$\sigma((\varepsilon_i)_{i\ge k+1}, (\zeta_z)_{(z-x)\cdot
e_1>0} )$-measurable for $y\in\mathcal{H}_1^x$.
Hence they are independent under $\mc P\otimes Q_\beta$ and we have
%
\begin{eqnarray}
\label{e7} && \bar{\mathbb P}_\lambda(B\circ\bar{\theta}_{\tau_1}
\cap A)
\nonumber\\[-6pt]\\[-10pt]\nonumber
&&\qquad = \sum_{k\ge1}\bar{\mathbb P}_\lambda
\bigl(A\cap\{\tilde{S}_k<\infty\}\bigr) E_P \biggl[
\sum_y \nu_{\omega^\lambda,1}(X_{T_1}=y)
\bar{P}_{\omega^\lambda}^y\bigl(B\cap\{T_{-1}=\infty\}\bigr)
\biggr].\hspace*{-15pt}
\nonumber
\end{eqnarray}
Substituting $B$ with the set of all events, we get
%
\begin{equation}
\label{e8} \bar{\mathbb P}_\lambda(A)= \sum
_{k\ge1}\bar{\mathbb P}_\lambda\bigl(A\cap\{
\tilde{S}_k<\infty\}\bigr) E_{\mc P} \biggl[ \sum
_y \mu_{\omega^\lambda,1}(y) \bar{P}_{\omega^\lambda}^y(T_{-1}=
\infty) \biggr].
\end{equation}
Equalities (\ref{e7}) and (\ref{e8}) yield
\[
\bar{\mathbb P}_\lambda(B\circ\bar{\theta}_{\tau_1}\mid A) =
\frac{E_{\mc P} [
\sum_y
\mu_{\omega^\lambda,1}(y)
\bar{P}_{\omega^\lambda}^y(B\cap\{T_{-1}=\infty\}) ]}{
E_{\mc P} [
\sum_y
\mu_{\omega^\lambda,1}(y)
\bar{P}_{\omega^\lambda}^y(T_{-1}=\infty) ]
}.
\]
The lemma is proved for the case $k=1$. The general case $k>1$ follows
by induction.
\end{pf}

We say that a sequence of random variables $(Y_i)_{i\in\mathbb{N}}$
is $m$-\textit{dependent} ($m\in\mathbb{N}$) if
\[
\sigma(Y_i; 1\le i\le n)\quad\mbox{and}\quad\sigma(Y_j; j>
n+m)\qquad\mbox{are independent }\forall n\in\mb N.
\]
The law of large numbers and central limit theorem also hold for a
stationary \mbox{$m$-}dependent sequence with finite means and variances, see
\cite{Bi56}, Theorem 5.2.
The following proposition is an immediate consequence of Lemma~\ref{l5}.

\begin{proposition}\label{prop1}
Under $\bar{\mathbb P}_\lambda$, $(X_{\tau_{n+1}}-X_{\tau_n})_{n\ge
1}$ and $(\tau_{n+1}-\tau_n)_{n\ge1}$ are stationary 1-dependent
sequences. Furthermore, for all $n\ge1$, $(X_{\tau_{n+1}}-X_{\tau
_n},\break  \tau_{n+1}-\tau_n)$ has law
\begin{eqnarray*}
&& \anp(X_{\tau_{n+1}}-X_{\tau_n}\in\cdot,\tau_{n+1}-
\tau_n\in\cdot)
\\
&&\qquad = E_{\mc P} \biggl[ \sum_y
\mu_{\omega^\lambda,1}(y) \bar{P}_{\omega^\lambda}^y(X_{\tau_1}\in
\cdot,\tau_1\in\cdot,T_{-1}=\infty) \biggr]\Big/p_\lambda.
\end{eqnarray*}
\end{proposition}

\begin{pf}
For $k\ge0$, let
\[
\mc F_k=\sigma(\tau_k, X_1,\ldots,
X_{\tau_k}).
\]
Then, for $n\ge1$, $\mc F_{n-1}\subset\mc G_n$ and
\begin{eqnarray*}
&&\anp(X_{\tau_{n+1}}-X_{\tau_n}\in\cdot,\tau_{n+1}-
\tau_n\in\cdot\mid\mc F_{n-1})
\\
&&\qquad =\ane\bigl[ \anp( X_{\tau_{n+1}}-X_{\tau_n}\in\cdot,
\tau_{n+1}-\tau_n\in\cdot\mid\mc G_n )\mid\mc
F_{n-1} \bigr].
\end{eqnarray*}
By Lemma~\ref{l5}, the proposition is proved.
\end{pf}


\subsection{Moment estimates}\label{secmo}
We will show that the typical values of $e_1\cdot(X_{\tau
_{k+1}}-X_{\tau_k})$ and $\tau_{k+1}-\tau_k$, ($k\ge0$) are $C/(\beta\lambda)$ and
$C/(\beta\lambda^2)$, respectively.

%
\begin{theorem}\label{regdist}
Assume \textup{(i)}. There exists a constant $c>0$ such that
\[
\bmb{E}_\lambda\bigl[\exp(c\beta\lambda X_{\tau_1}) \bigr]<C
\]
for all $\lambda\in(0,\lambda_0)$ and $\beta\in(0,1)$.
\end{theorem}

\begin{pf}
Our proof contains several steps.
\begin{longlist}[2.]
\item[1.]
For $0\le k\le K-1$, set
\[
L_{k+1}=\inf\{n\ge\lambda_1 M_k\dvtx
\varepsilon_n=1\}-\lambda_1 M_k+1.
\]
Then $L_1$ is the number of coins tossed to get the first ``$1$'' and
%
\begin{equation}
\label{e33} X_{S_1}\cdot e_1=L_1/
\lambda_1.
\end{equation}
Since $(L_i)_{i\ge1}$ depends only on the coins $(\varepsilon_i)_{i\ge0}$,
it is easily seen that they are i.i.d. geometric with
parameter $\beta$. Hence, for $i\ge1$, $s\in(0,1)$,
%
\begin{equation}
\label{e11} \bar{E}_{\omega^\lambda} \bigl[e^{s\beta L_i}\bigr] =
\frac{\beta e^{s\beta}}{1-(1-\beta)e^{s\beta}}\le\frac{1}{1-s}.
\end{equation}
Moreover,
for $1\le k\le K-1$,
\[
(X_{S_{k+1}}-X_{S_k})\cdot e_1=N\circ
\theta_{S_k}+L_{k+1}/\lambda_1.
\]
\item[2.]
We will show that
%
\begin{equation}
\label{e10} \bmb{E}_\lambda\bigl[\exp(s\beta\lambda_1
X_{\tau_1}\cdot e_1) \bigr] =\sum
_{k\ge1}\ane\bigl[\exp(s\beta\lambda_1
X_{S_k}\cdot e_1), S_k<\infty
\bigr]p_\lambda.
\end{equation}
By the definition of $X_{\tau_1}$,
\begin{eqnarray*}
&&\bmb{E}_\lambda\bigl[\exp(s\beta\lambda_1
X_{\tau_1}\cdot e_1) \bigr]
\\
&&\qquad =\sum_{k\ge1,x,y}E_{\mc P} \bigl[
\bar{E}_{\omega^\lambda} \bigl[\exp(s\beta\lambda_1 X_{\tilde
{S}_k}
\cdot e_1+s\beta),\tilde S_k<\infty, X_{\tilde S_k}=x
\bigr]
\\
&&\hspace*{149pt}{} \times\mu_{\omega^\lambda,1}^x(x+y)P_{\omega^\lambda
}^{x+y}(T_{-1}=
\infty) \bigr].
\end{eqnarray*}
Since
\[
\bar{E}_{\omega^\lambda} \bigl[\exp(s\beta\lambda_1 X_{\tilde
{S}_k}
\cdot e_1+s\beta), \tilde S_k<\infty, X_{\tilde S_k}=x
\bigr]
\]
is $\sigma(\zeta_z\dvtx  z\cdot e_1<x\cdot e_1)$-measurable, and $\mu
_{\omega^\lambda,1}^x(x+y)P_{\omega^\lambda}^{x+y}(T_{-1}=\infty)$
is $\sigma(\zeta_z\dvtx z\cdot e_1\ge x\cdot e_1)$-measurable, they are
independent under $\mc P$.
Therefore,
%
\begin{eqnarray}
\label{e22} &&\ane\bigl[\exp(s\beta\lambda_1 X_{\tau_1}
\cdot e_1) \bigr]
\nonumber
\\
&&\qquad =\sum_{k\ge1, x,y}\ane\bigl[\exp(s\beta
\lambda_1 X_{\tilde
{S}_k}\cdot e_1+s\beta), \tilde
S_k<\infty, X_{\tilde S_k}=x \bigr]
\\
&&\hspace*{60pt}{} \times E_{\mc P} \bigl[\mu_{\omega^\lambda,1}^x(x+y)P_{\omega
^\lambda}^{x+y}(T_{-1}=
\infty) \bigr].\nonumber
\end{eqnarray}
Equation (\ref{e10}) follows.
\item[3.]
Next, we will show that for $k\ge1$,
%
\begin{eqnarray}
\label{e12} && \ane\bigl[\exp(s\beta\lambda_1 X_{S_k}
\cdot e_1), S_k<\infty\bigr]
\nonumber\\[-8pt]\\[-8pt]\nonumber
&&\qquad \le\biggl(\frac{A(s,\lambda,\beta)}{1-s} \biggr)^{k-1}\ane\bigl[\exp
(s\beta
\lambda_1 X_{S_1}\cdot e_1), S_1<
\infty\bigr],\nonumber
\end{eqnarray}
where
\[
A(s,\lambda,\beta):=E_{\mc P} \biggl[\sum_y
\mu_{\omega^\lambda
,1}(y)E_{\omega^\lambda}^y\bigl[e^{s\beta N},
T_{-1}<\infty\bigr] \biggr].
\]
By definition,
\begin{eqnarray*}
&&\ane\bigl[\exp(s\beta\lambda_1 X_{S_{k+1}}\cdot
e_1), S_{k+1}<\infty\bigr]
\\
&&\qquad =\ane\bigl[ \exp(s\beta\lambda_1 X_{S_k}\cdot
e_1+s\beta N\circ\theta_{S_k}+s\beta L_{k+1}),S_k<
\infty,T_{-1}\circ\theta_{S_k}<\infty\bigr].
\end{eqnarray*}
Noting that $L_{k+1}$ is independent of $\sigma\{R_k, X_1,\ldots,
X_{R_k}\}$, we get
%
\begin{eqnarray}
\label{e26} && \ane\bigl[\exp(s\beta\lambda_1 X_{S_{k+1}}
\cdot e_1), S_{k+1}<\infty\bigr]\nonumber
\\
&&\qquad =\bmb{E}_\lambda\bigl[ \exp(s\beta\lambda_1
X_{S_k}\cdot e_1+s\beta N\circ\theta_{S_k}),S_k<
\infty,T_{-1}\circ\theta_{S_k}<\infty\bigr]
\nonumber\\[-8pt]\\[-8pt]\nonumber
&&\quad\qquad{}\times \ane\bigl[e^{s\beta L_{k+1}}\bigr]
\\
&&\hspace*{-3pt}\qquad \stackrel{\fontsize{8.36pt}{8.36pt}{(\ref{e11})}} {\le} \ane\bigl[ \exp(s\beta
\lambda_1 X_{S_k}\cdot e_1+s\beta N\circ\theta
_{S_k}),S_k<\infty,T_{-1}\circ
\theta_{S_k}<\infty\bigr]/(1-s).\hspace*{-5pt}
\nonumber
\end{eqnarray}
Further, by the same argument as in (\ref{e22}),
\begin{eqnarray*}
&&\bmb{E}_\lambda\bigl[ \exp(s\beta\lambda_1
X_{S_k}\cdot e_1+s\beta N\circ\theta_{S_k}),S_k<
\infty,T_{-1}\circ\theta_{S_k}<\infty\bigr]
\\
&&\qquad =\sum_{x,y}\bmb{E}_\lambda\bigl[\exp(s
\beta\lambda_1 x\cdot e_1+s\beta),\tilde S_k<
\infty, X_{\tilde S_k}=x \bigr]
\\
&&\hspace*{45pt}{}\times E_{\mc P} \bigl[\mu_{\omega^\lambda,
1}^x(x+y)E_{\omega^\lambda}^{x+y}
\bigl[e^{s\beta
N},T_{-1}<\infty\bigr] \bigr]
\\
&&\qquad =\bmb{E}_\lambda\bigl[\exp(s\beta\lambda_1
X_{S_k}\cdot e_1), S_k<\infty\bigr]A(s,
\lambda,\beta).
\end{eqnarray*}
Combining the above equality and (\ref{e26}), ineqaulity (\ref{e12})
follows by induction.
\item[4.]
By (\ref{e10}) and (\ref{e12}), we have
\[
\bmb{E}_\lambda\bigl[\exp(s\beta\lambda_1 X_{\tau_1}
\cdot e_1) \bigr] \le p_\lambda\ane\bigl[\exp(s\beta
\lambda_1 X_{S_1}\cdot e_1), S_1<
\infty\bigr] \sum_{k=0}^\infty\biggl(
\frac{A(s,\lambda,\beta)}{1-s} \biggr)^k.
\]
Since $p_\lambda\le1$, and
[by (\ref{e33}) and (\ref{e11})]
\[
\bmb{E}_\lambda\bigl[\exp(s\beta\lambda_1 X_{S_1}
\cdot e_1),S_1<\infty\bigr] =\ane\bigl[e^{s\beta L_1}
\bigr] \le\frac{1}{1-s},
\]
to prove Theorem~\ref{regdist}, we only need that when $s>0$ is small enough,
%
\begin{equation}
\label{e31} A(s,\lambda,\beta)<C<1.
\end{equation}
For any $m\in\mb N$,
\begin{eqnarray*}
A(s,\lambda,\beta) &\le& e^{s\beta m}\mb P_\lambda(T_{-1}<
\infty)
\\
&&{}  +\sum_{n=m}^\infty
e^{s\beta n}E_{\mc P} \biggl[\sum_y
\mu_{\omega
^\lambda,1}(y)P_{\omega^\lambda}^y(N=n, T_{-1}<
\infty) \biggr].
\end{eqnarray*}
Hence [note that $\mb P_\lambda(T_{-1}<\infty)\le e^{-\rho/2}$], to
prove (\ref{e31}), it suffices to show
%
\begin{equation}
\label{e32} E_{\mc P} \biggl[\sum_y
\mu_{\omega^\lambda,1}(y)P_{\omega^\lambda
}^y(N=n, T_{-1}<
\infty) \biggr]<Ce^{-cn}.
\end{equation}
\item[5.] Recall the definition of $N$ in (\ref{e46}):
\begin{eqnarray*}
&&E_{\mc P} \biggl[\sum_y
\mu_{\omega^\lambda,1}(y)P_{\omega^\lambda
}^y(N=n, T_{-1}<
\infty) \biggr]
\\
&&\qquad \le CE_{\mc P} \biggl[\sum_yP_{\omega^\lambda}(X_{T_1}=y)P_{\omega
^\lambda}^y(N=n,
T_{-1}<\infty) \biggr]
\\
&&\qquad \le CE_{\mc P} \biggl[\sum_{y,z}P_{\omega^\lambda
}(X_{T_1}=y)P_{\omega^\lambda}^y(X_{T_{n-1}}=z)P_{\omega^\lambda
}^z(T_{-n-1}<T_1)
\biggr]
\\
&&\qquad =CE_{\mc P} \biggl[\sum_{z}P_{\omega^\lambda}(X_{T_n}=z)P_{\omega
^\lambda}^z(T_{-n-1}<T_1)
\biggr].
\end{eqnarray*}
For $k\ge0$, let $z_k:=(ne_1+ke_2)/\lambda_1$ and
\[
A_k:= \biggl\{z\in\mc H\dvtx \frac{k}{\lambda_1}\le\biggl\llvert
X_{T_n}-\frac
{n}{\lambda_1}e_1\biggr\rrvert<
\frac{k+1}{\lambda_1}\biggr\}.
\]
Then by the Harnack inequality, for any $x\in A_k$,
\[
P_{\omega^\lambda}^x(T_{-n-1}<T_1)\le
C_2 P_{\omega^\lambda
}^{z_k}(T_{-n-1}<T_1),
\]
where $C_2$ is the constant in (\ref{e49}). Hence,
\begin{eqnarray*}
&&E_{\mc P} \biggl[\sum_{z}P_{\omega^\lambda}(X_{T_n}=z)P_{\omega
^\lambda}^z(T_{-n-1}<T_1)
\biggr]
\\
&&\qquad \le\mb P_\lambda\biggl(\biggl\llvert X_{T_n}-
\frac{n}{\lambda
_1}e_1\biggr\rrvert\ge\frac{n}{\lambda_1} \biggr)
\\
&&\quad\qquad{} +C_2\sum_{k=0}^{n-1}E_{\mc P}
\biggl[\sum_{z\in A_k}P_{\omega^\lambda
}(X_{T_n}=z)P_{\omega^\lambda}^{z_k}(T_{-n-1}<T_1)
\biggr]
\\
&&\hspace*{-14pt}\qquad \stackrel{\fontsize{8.36pt}{8.36pt}{\mathrm{Lemma}~\ref{prop2}}} {\le} Ce^{-cn}
+C_2E_{\mc P}
\Biggl[\sum_{k=0}^{n-1}P_{\omega^\lambda}^{z_k}(T_{-n-1}<T_1)
\Biggr]
\\
&&\qquad =Ce^{-cn}+C_2n\mb{P}_\lambda(T_{-n-1}<T_1)
\le Ce^{-cn}.
\end{eqnarray*}
Inequality (\ref{e32}) is proved.\quad\qed
\end{longlist}\noqed
\end{pf}

%
\begin{corollary}\label{cor1}
Assume \textup{(i)}.
For all $n\ge0$, $\lambda\in(0,\lambda_0)$ and $\beta\in(0,1)$,
%
\begin{eqnarray}\label{e13}
\ane\bigl[\exp\bigl(c_4\beta\lambda(X_{\tau_{n+1}}-X_{\tau
_n})\cdot e_1 \bigr) \bigr]&<&C,
\\
\label{e16} \anp\bigl( \beta\lambda_1^2(
\tau_{n+1}-\tau_n)>t \bigr)&\le& C e^{-c\sqrt t}\qquad
\forall t>0.
\end{eqnarray}
Here, $c_4>0$ is a constant.
\end{corollary}

\begin{pf}
1.~First, we consider the case $n=0$, $\tau_{n+1}-\tau_n=\tau_1$.
Ineqaulity (\ref{e13}) is the conclusion of Theorem~\ref{regdist}.
To prove (\ref{e16}), note that for any $m\in\mb N$,
%
\begin{eqnarray}
\label{e34} &&\bmb{P}_\lambda\bigl(\beta\lambda_1^2
\tau_1>t\bigr) \le\bmb{P}_\lambda(X_{\tau_1}\cdot
e_1\ge m/\lambda_1)+\bmb{P}_\lambda\bigl(\beta
\lambda_1^2T_m>t\bigr).
\end{eqnarray}
By Theorem~\ref{regdist},
\[
\bmb{P}_\lambda(X_{\tau_1}\cdot e_1\ge m/
\lambda_1) \le Ce^{-c\beta m}.
\]
By Proposition~\ref{prop5} and Theorem~\ref{44444},
\begin{eqnarray*}
\anp\bigl(\beta\lambda_1^2T_m>t\bigr) &\le&
\anp(T_m>T_{-m})+\anp\bigl(\beta\lambda_1^2
\tilde T_m>t\bigr)
\\
&\le& Ce^{-cm}+Ce^{-ct/(\beta m)}.
\end{eqnarray*}
Coming back to (\ref{e34}), we get
\[
\bmb{P}_\lambda\bigl(\beta\lambda_1^2
\tau_1>t\bigr) \le Ce^{-c\beta m}+Ce^{-ct/(\beta m)}\qquad\forall m
\in\mb N.
\]
Ineqaulity (\ref{e16}) (for $n=0$) follows by letting $m=\lfloor\sqrt
t/\beta\rfloor$.

2.
Next, we will prove (\ref{e13}) for $n\ge1$.
By Lemma~\ref{l5},
\begin{eqnarray*}
&&\ane\bigl[\exp\bigl(s\beta\lambda_1 (X_{\tau_{n+1}}-X_{\tau_n})
\cdot e_1 \bigr) \bigr]
\\
&&\qquad = E_{\mc P} \biggl[ \sum_y
\mu_{\omega^\lambda,1}(y)\bar{E}_{\omega^\lambda}^y\bigl[\exp(s\beta
\lambda_1 X_{\tau_1}\cdot e_1 ),T_{-1}=
\infty\bigr] \biggr]\Big/p_\lambda
\\
&&\qquad \le E_{\mc P} \biggl[ \sum_y
\mu_{\omega^\lambda,1}(y)\bar{E}_{\omega^\lambda}^y\bigl[\exp(s\beta
\lambda_1 X_{\tau_1}\cdot e_1 )\bigr]
\biggr]\Big/p_\lambda.
\end{eqnarray*}
By the same argument as in (\ref{e10}) and (\ref{e12}), we get
\begin{eqnarray*}
&& E_{\mc P} \biggl[ \sum_y
\mu_{\omega^\lambda,1}(y)\bar{E}_{\omega^\lambda}^y\bigl[\exp(s\beta
\lambda_1 X_{\tau_1}\cdot e_1 )\bigr] \biggr]
\\
&&\qquad =\sum_{k\ge1}E_{\mc P} \biggl[\sum
_y\mu_{\omega^\lambda,1}(y)\bar{E}_{\omega^\lambda}^y
\bigl[\exp(s\beta\lambda_1 X_{S_k}\cdot e_1),S_k<
\infty\bigr] \biggr]p_\lambda,
\end{eqnarray*}
and
\begin{eqnarray*}
&& E_{\mc P} \biggl[\sum_y
\mu_{\omega^\lambda,1}(y)\bar{E}_{\omega
^\lambda}^y\bigl[\exp(s\beta
\lambda_1 X_{S_k}\cdot e_1),S_k<
\infty\bigr] \biggr]
\\
&&\qquad \le\biggl(\frac{A(s,\lambda,\beta)}{1-s} \biggr)^{k-1}E_{\mc P} \biggl[
\sum_y\mu_{\omega^\lambda,1}(y)\bar{E}_{\omega^\lambda}^y
\bigl[\exp(s\beta\lambda_1 X_{S_1}\cdot
e_1),S_1<\infty\bigr] \biggr]
\\
&&\hspace*{-13pt}\qquad \stackrel{\fontsize{8.36pt}{8.36pt}{(\ref{e33}),\,(\ref{e11})}}{\le} \frac{A(s,\lambda,\beta
)^{k-1}}{(1-s)^k}\qquad
\forall k\ge1.
\end{eqnarray*}
Therefore,
\[
\ane\bigl[\exp\bigl(s\beta\lambda_1 (X_{\tau_{n+1}}-X_{\tau_n})
\cdot e_1 \bigr) \bigr] \le\sum_{k\ge1}
\frac{A(s,\lambda,\beta)^{k-1}}{(1-s)^k}.
\]
By (\ref{e32}), ineqaulity (\ref{e13}) is proved.

3. Finally, we will prove (\ref{e16}) for $n\ge1$. Similar to (\ref
{e34}), for any $m\in\mb N$,
%
\begin{eqnarray}
\label{e35} &&\anp\bigl(\beta\lambda_1^2(
\tau_{n+1}-\tau_n)>t\bigr)
\nonumber\\[-8pt]\\[-8pt]\nonumber
&&\qquad \le\anp\bigl((X_{\tau_{n+1}}-X_{\tau_n})\cdot e_1\ge m/
\lambda_1\bigr)+\anp\bigl(\beta\lambda_1^2
T_m\circ\theta_{\tau_n}>t\bigr).
\nonumber
\end{eqnarray}
By (\ref{e13}),
%
\begin{equation}
\label{e36} \anp\bigl((X_{\tau_{n+1}}-X_{\tau_n})\cdot
e_1\ge m/\lambda_1\bigr) \le C e^{-cm\beta}\qquad
\forall m\in\mb N.
\end{equation}
By Lemma~\ref{l5},
\begin{eqnarray*}
&&\anp\bigl(\beta\lambda_1^2 T_{m}\circ
\theta_{\tau_n}>t\bigr)
\\
&&\qquad = \sum_y E_{\mc P} \bigl[
\mu_{\omega^\lambda,1}(y) \bar{P}_{\omega^\lambda}^y\bigl(\beta
\lambda_1^2 T_{m}>t, T_{-1}=\infty
\bigr) \bigr]/p_\lambda
\\
&&\qquad \le\sum_y E_{\mc P} \bigl[
\mu_{\omega^\lambda,1}(y) \bar{P}_{\omega^\lambda}^y\bigl(\beta
\lambda_1^2 \tilde T_{m}>t\bigr)
\bigr]/p_\lambda.
\end{eqnarray*}
Applying Theorem~\ref{44444} to the above inequality, we have
%
\begin{equation}
\label{e37} \anp\bigl(\beta\lambda_1^2
T_{m}\circ\theta_{\tau_n}>t\bigr) \le C e^{-ct/(m\beta)}/p_\lambda
\qquad\forall m\in\mb N.
\end{equation}
Combining (\ref{e35}), (\ref{e36}) and (\ref{e37}) and letting
$m=\lfloor\sqrt t/\beta\rfloor$,
\[
\anp\bigl(\beta\lambda_1^2(\tau_{n+1}-
\tau_n)>t\bigr) \le C e^{-c\sqrt t}/p_\lambda.
\]
It remains to show that
\[
p_\lambda>C>0\qquad\forall\lambda\in(0,\lambda).
\]
By (\ref{e17}) and (\ref{e18}),
\begin{eqnarray*}
p_\lambda&=& E_{\mc P} \bigl[P_{\omega^\lambda}^{0.5e_1/\lambda_1}
(T_{-0.5}=\infty\mid T_{0.5}<T_{-0.5} ) \bigr]
\\
&=&E_{\mc P} \bigl[P_{\omega^\lambda} (T_{-0.5}=\infty\mid
T_{0.5}<T_{-0.5} ) \bigr]
\\
&\ge&\mb{P}_\lambda(T_{-0.5}=\infty)\stackrel{\fontsize{8.36pt}{8.36pt}{\mathrm{Proposition~\ref{prop5}}}} {\ge} C\bigl(1-e^{-\rho/4}\bigr).
\end{eqnarray*}
%

Our proof of (\ref{e16}) is complete.
\end{pf}

By Corollary \ref{cor1},
we conclude that for any $p\ge1, k\ge0$, there exists a constant
$C(p)<\infty$ such that
%
\begin{eqnarray}
\bar{\mathbb E}_\lambda\bigl[\bigl(\beta\lambda_1^2
\tau_1\circ\theta_{\tau_k}\bigr)^p\bigr]
&<&C(p)\label{e14}
\end{eqnarray}
and
\begin{eqnarray}
\bar{\mathbb E}_\lambda\bigl[(\beta
\lambda_1 X_{\tau_1}\circ\theta_{\tau_k})^p
\bigr] &<&C(p).\label{e15}
\end{eqnarray}
Moreover,
by the law of large numbers,
\[
v_\lambda\stackrel{\bar{\mathbb P}_\lambda{\mbox{\fontsize{8.36pt}{8.36pt}{-a.s.}}}} {=}
\lim_{n\to\infty}\frac{X_{\tau_n}\cdot e_1}{\tau_n}=\frac{\ane
[X_{\tau_2}-X_{\tau_1}]}{\ane[\tau_2-\tau_1]}.
\]
On the other hand,
\begin{eqnarray*}
\ane\bigl[(X_{\tau_2}-X_{\tau_1})\cdot e_1\bigr] &\ge&
1/\lambda_1,
\\
v_\lambda\cdot e_1&=&\lambda E_{\mc Q_\lambda}\bigl[d(\xi)
\bigr]\cdot e_1\le2\lambda.
\end{eqnarray*}
Hence,
%
\begin{equation}
\label{e29} \ane[\tau_2-\tau_1]\ge\frac{1/\lambda_1}{2\lambda}
\ge C/\lambda^2.
\end{equation}
%
\section{Proof of Theorem~\texorpdfstring{\protect\ref{ER0}}{1}}\label{secpro}
Let $\alpha_n=\alpha_n(\beta,\lambda):=\bar{\mathbb E}_\lambda
\tau_n$.
Note that by (\ref{e14}) and~(\ref{e29}),
%
\begin{equation}
\label{e20} \frac{Cn}{\lambda^2} \le\alpha_n \le\frac{C'n}{\lambda^2}.
\end{equation}
%

\begin{lemma}\label{l6}
Assume \textup{(i)}. Let $f$ be a function that satisfies (\ref{f}). Then for
$\beta>0$ and $\lambda\in(0,1/N_f)$,
\[
\Biggl\llvert E_{\mc Q_\lambda} f-\frac{1}{\bmb{E}_\lambda\tau_n}\ane
\Biggl[\sum
_{i=0}^{\tau_n} f(\bar\zeta_i) \Biggr] \Biggr
\rrvert\le C\llVert f\rrVert_\infty/n\qquad\mbox{for all }n\in\mb N.
\]
\end{lemma}

\begin{pf}
The lemma is trivial when $n=1$, so we only consider $n\ge2$.
Recall that $\tau_0=0$. For $k\ge0$, set
\[
Z_k=Z_k(f) =\sum_{i=\tau_k}^{\tau_{k+1}-1}f(
\bar\zeta_i).
\]
Since $N_f\le\frac{1}{\lambda}\le\ane[X_{\tau_2}-X_{\tau_1}]$,
we see that
$(Z_k)_{k\ge0}$ is an 1-dependent sequence.
On one hand,
%
\begin{equation}
\label{e28} \frac{1}{\bmb{E}_\lambda\tau_n}\bmb{E}_\lambda\Biggl[\sum
_{i=0}^{\tau_n} f(\bar\zeta_i) \Biggr] =
\frac{(n-1)\ane Z_1+\ane Z_0}{(n-1)\ane[\tau_2-\tau_1]+\ane\tau_1}.
\end{equation}
On the other hand, since the $\mb{P}_\lambda$-law of $\bar\zeta_n$
converges weakly to $\mc Q_\lambda$, by (\ref{e47}),
\[
\mc Q_\lambda f=\lim_{n\to\infty} \bmb{E}_\lambda
\Biggl[ \frac{1}{n}\sum_{i=0}^{n-1}f(
\bar\zeta_i) \Biggr].
\]
Hence, by the law of large numbers,
%
\begin{equation}
\label{e21} \mc Q_\lambda f=\frac{\ane[Z_1]}{\ane[\tau_2-\tau_1]}.
\end{equation}
%
The lemma follows by combining (\ref{e28}), (\ref{e21}) and using the
moment bounds (\ref{e14}) and (\ref{e29}).
\end{pf}


\begin{lemma}\label{thm3}
Assume \textup{(i)}. Let $f$ be a function that satisfies (\ref{f}).
Then
\[
\frac{1}{\alpha_n} \bmb{E}_\lambda\Biggl\llvert\sum
_{i=0}^{\tau_n}f(\bar\zeta_i) -\sum
_{i=0}^{\alpha_n}f(\bar\zeta_i) \Biggr\rrvert
\le\frac{C\llVert f\rrVert_\infty}{\sqrt n} \qquad\forall n\in\mathbb
{N}, \lambda\in(0,1/N_f).
\]
\end{lemma}

\begin{pf}
Noting that the left-hand side is less than
\[
\frac{\llVert f\rrVert_\infty}{\alpha_n}\bmb{E}_\lambda\llvert\tau_n-
\alpha_n\rrvert\le\frac{\llVert f\rrVert_\infty}{\alpha_n}\sqrt{\var
\tau_n},
\]
by (\ref{e20}), it suffices to show
%
\begin{equation}
\label{e23} \var\tau_n\le Cn/\lambda^4.
\end{equation}
Indeed, by the inequality $(a+b)^2\le2(a^2+b^2)$, we have
\begin{eqnarray*}
\var\tau_n &=& \var\Biggl[\sum_{k=0}^{n-1}(
\tau_k-\tau_{k-1})\Biggr]
\\
&\le& 2 \Biggl(\var\Biggl[\sum_{k=0}^{\lfloor (n-1)/2\rfloor}(
\tau_{2k+1}-\tau_{2k}) \Biggr] +\var\Biggl[\sum
_{k=1}^{\lfloor n/2\rfloor}(\tau_{2k}-\tau
_{2k-1}) \Biggr] \Biggr).
\end{eqnarray*}
Since $(\tau_{2k+1}-\tau_{2k})_{k\ge0}$ and $(\tau_{2k}-\tau
_{2k-1})_{k\ge1}$ are i.i.d. sequences,
we conclude that
\[
\var\tau_n \le2\sum_{k=0}^{\lfloor (n-1)/2\rfloor}
\var[\tau_{2k+1}-\tau_{2k}] +2\sum
_{k=1}^{\lfloor n/2\rfloor}\var[\tau_{2k}-
\tau_{2k-1}] \stackrel{\fontsize{8.36pt}{8.36pt}{(\ref{e14})}} {\le} Cn/\lambda^4.
\]
This completes the proof of (\ref{e23}).
\end{pf}

\begin{pf*}{Proof of Theorem~\ref{ER2}}
Since the left-hand side
is uniformly bounded (by $2\llVert f\rrVert_\infty$) for all
$t$, the case $t<\lambda^2\alpha_1\le C$ is trivial.
For $t\ge\lambda^2\alpha_1$, we let $n=n(t,\lambda)\ge1$ be the
integer that satisfies
\[
\alpha_n\le\frac{t}{\lambda^2}<\alpha_{n+1}.
\]

Since
\begin{eqnarray*}
&&\Biggl\llvert\frac{\lambda^2}{t}\ane\sum_{i=0}^{\lceil t/\lambda
^2\rceil
}f(
\bar\zeta_i) -\frac{1}{\alpha_n}\ane\sum_{i=0}^{\alpha_n}f(
\bar\zeta_i) \Biggr\rrvert
\\
&&\qquad \le\Biggl\llvert\frac{\lambda^2}{t}\ane\Biggl[\sum
_{i=0}^{\lceil t/\lambda^2\rceil}f(\bar\zeta_i) -\sum
_{i=0}^{\alpha_n}f(\bar\zeta_i) \Biggr] \Biggr
\rrvert+\Biggl\llvert\biggl(\frac{\lambda^2}{t}-\frac{1}{\alpha
_n}\biggr)\ane\sum
_{i=0}^{\alpha_n}f(\bar\zeta_i)\Biggr
\rrvert
\\
&&\qquad \le\llVert f\rrVert_\infty\biggl(\frac{\lambda^2}{t}(\alpha
_{n+1}-\alpha_n)+\biggl(\frac{1}{\alpha_n}-
\frac{1}{\alpha_{n+1}}\biggr)\alpha_n \biggr)
\\
&&\hspace*{-3.6pt}\qquad \stackrel{\fontsize{8.36pt}{8.36pt}{(\ref{e20})}} {\le} C\llVert f\rrVert_\infty/n.
\end{eqnarray*}
Theorem~\ref{ER2} follows by combining Lemma \ref{l6}, Lemma \ref
{thm3} and the above inequality.
\end{pf*}

\begin{pf*}{Proof of Theorem~\ref{ER0}}
Since the space $\Omega$ of the environment is compact under the
product topology, it suffices to show that
\[
\lim_{\lambda\to0}\mc Q_\lambda f=\mc Q f
\]
for all f that satisfies (\ref{f}). The above equality follows
immediately from Theorems~\ref{ER1}~and~\ref{ER2}.
\end{pf*}

\section{Proof of Theorem~\texorpdfstring{\protect\ref{thm6}}{2}}\label{secball}
Let us recall the regeneration structure defined by Sznitman and Zerner
\cite{SZ}.
For a path $(X_n)_{n\ge0}$, we call $t>0$ a \textit{renewal time}\footnote{It is usually called a \textit{regeneration time} in the RWRE
literature. But we use a different name to distinguish with the
regeneration structure defined in Section~\ref{secreg}.} in the
direction $\ell$ if
\[
X_m\cdot\ell<X_t\cdot\ell<X_n\cdot\ell
\]
for all $m,n$ such that $m<t<n$.
For ballistic RWRE, the renewal times exist a.s. and have finite first
moments. We let
\[
T(1)<T(2)<\cdots
\]
denote all the renewal times. Then $(X_{T(k+1)}-X_{T(k)},
T(k+1)-T(k))_{k\ge1}$ is an i.i.d. sequence under $\mb P$.


%
\begin{lemma}\label{l3}
If the $\mc P$-law of $\omega$ satisfies Sznitman's \textup{(T$'$)} condition,
then there exists a constant $\lambda_0>0$ such that for all $\lambda
\in[0,\lambda_0)$:
\begin{longlist}[a]
\item[(a)] $\omega^\lambda$ satisfies \textup{(T$'$)};
\item[(b)] $\anne[T(1)^2]<C$ and $\anne[(T(2)-T(1))^2]<C$.
\end{longlist}
\end{lemma}

\begin{pf}
It is shown in \cite{BDR}, Theorem 1.6, that \textup{(T$'$)} is equivalent to a
polynomial ballisticity condition ($\mc P$). Note that ($\mc P$) only
involves checking a strict inequality for some (finitely many) exit
probabilities from a finite box (see \cite{BDR}, Definition~1.4).
Hence, there exists $\lambda_0>0$ such that ($\mc P$) holds for all
$\omega^\lambda$, $\lambda\in[0,\lambda_0)$, with the same
constants in the upper bounds of \cite{BDR}, Definition 3.2. We have
proved~(a). Furthermore, by \cite{Sz3}, Proposition 3.1 and \cite{BDR},
Theorem 1.6, $(\mc P)$ implies that the regeneration time has
finite moments. Therefore, the second moments of $T(1)$ and $T(2)-T(1)$
(under $\annp$) can be bounded by the same constant [since they are
deduced from the same ($\mc P$) condition] for all $\lambda\in
[0,\lambda_0)$. (b) is proved.
\end{pf}

%
\begin{theorem}\label{thm5}
Assume that $\omega$ satisfies Sznitman's \textup{(T$'$)} condition. If $f$
satisfies~(\ref{f}), then for any $t\ge1$, $n\in\mb N$ and $\lambda
\in[0,\lambda_0)$,
\[
\anne\Biggl[ \Biggl(\sum_{i=0}^n
\bigl(f(\bar\zeta_i)-\mathcal{Q}_\lambda f\bigr)
\Biggr)^2 \Biggr] \le CN_f\llVert f\rrVert
_\infty^2 n.
\]
\end{theorem}

To prove this theorem, we need two lemmas.

\begin{lemma}\label{l7}
Assume that $\omega$ satisfies Sznitman's \textup{(T$'$)} condition. If $f$
satisfies~(\ref{f}), then for any $\lambda\in[0,\lambda_0)$,
\[
\anne\Biggl[ \Biggl[ \sum_{i=0}^{T(n)}
\bigl(f(\bar\zeta_i)-\mathcal{Q}_\lambda f \bigr)
\Biggr]^2 \Biggr]\le CN_f\llVert f\rrVert
_\infty^2 n.
\]
\end{lemma}

\begin{pf}
For $k\ge0$, let
\[
Z_k=Z_k(f):=\sum_{i=T(k)}^{T(k+1)-1}
\bigl(f(\bar\zeta_i)-\mathcal{Q}_\lambda f \bigr).
\]
Then $(Z_k)_{k\ge N_f}$ is a $N_f$-dependent and stationary sequence.
Moreover, for $k\ge N_f$,
\[
\anne Z_k=0.
\]
Hence, for $n>N_f$,
\begin{eqnarray*}
&& \anne\Biggl[\Biggl(\sum_{k=N_f}^{n-1}
Z_{k}\Biggr)^2 \Biggr]
\\
&&\qquad =\sum_{k=N_f}^n E_\lambda
\bigl[Z_{k}^2\bigr]+2\sum_{j=N_f}^{n-N_f}
\sum_{k=j+1}^{N_f}E_\lambda[Z_j
Z_k]
\\
&&\qquad \le3n N_f E_\lambda\bigl[Z_{N_f}^2
\bigr]\le CnN_f \llVert f\rrVert_\infty^2.
\end{eqnarray*}
Noting that
\[
\anne\Biggl[\Biggl(\sum_{k=0}^{N_f-1}
Z_k\Biggr)^2 \Biggr] \le\llVert f\rrVert
_\infty^2\anne\bigl[T(N_f)^2\bigr]
\le CN_f\llVert f\rrVert_\infty^2,
\]
our proof is complete.
\end{pf}

%
\begin{lemma}\label{l8}
Let $\alpha_n=\alpha(n,\lambda)=\anne T(n)$.
Assume that $\omega$ satisfies Sznitman's \textup{(T$'$)} condition. If $f$
satisfies (\ref{f}), then for any $\lambda\in[0,\lambda_0')$,
\[
\anne\Biggl[ \Biggl( \sum_{i=0}^{T(n)}f(
\bar\zeta_i)- \sum_{i=0}^{\alpha_n}f(
\bar\zeta_i) \Biggr)^2 \Biggr] \le C\llVert f\rrVert
_\infty^2n.
\]
\end{lemma}

\begin{pf}
\begin{eqnarray*}
\anne\Biggl[ \Biggl(\sum_{i=0}^{T(n)}f(\bar
\zeta_i)- \sum_{i=0}^{\alpha_n}f(\bar
\zeta_i) \Biggr)^2 \Biggr] &\le&\llVert f\rrVert
_\infty^2\anne\bigl[\bigl(T(n)-\alpha_n
\bigr)^2\bigr]
\\
&\le&\llVert f\rrVert_\infty^2\sum
_{i=0}^{n-1}\var\bigl[T(i+1)-T(i)\bigr].
\end{eqnarray*}
By Lemma~\ref{l3}(b), the lemma follows.
\end{pf}

\begin{pf*}{Proof of Theorem~\ref{thm5}}
Set
\[
\tilde f(\zeta):=f(\zeta)-\mc Q f.
\]
By Lemmas~\ref{l7}~and~\ref{l8}, for any $m\in\mb N$,
\[
\anne\Biggl[ \Biggl(\sum_{i=0}^{\alpha_m}\tilde
f(\bar\zeta_i) \Biggr)^2 \Biggr] \le CN_f
\llVert f\rrVert_\infty^2m.
\]
For $n\ge1$, we let $m=m(s,\lambda)\ge0$ be the integer that satisfies
\[
\alpha_m\le n<\alpha_{m+1}.
\]
Thus,
\begin{eqnarray*}
&& \anne\Biggl[ \Biggl(\sum_{i=0}^{n}
\tilde f(\bar\zeta_i) \Biggr)^2 \Biggr]
\\
&&\qquad \le2\anne\Biggl[ \Biggl(\sum_{i=0}^{\alpha_m}
\tilde f(\bar\zeta_i) \Biggr)^2 \Biggr] +8\llVert f
\rrVert_\infty^2 (\alpha_{n+1}-
\alpha_n) 
\le C N_f\llVert f
\rrVert_\infty^2 n.
\end{eqnarray*}\upqed
\end{pf*}


\begin{pf*}{Proof of Theorem~\ref{thm6}}
Recall the definitions of $G(\cdot,\cdot)$ and $a(\zeta,e)$ in
Section~\ref{sec1}.
Since
\[
\bigl((A_n, B_n) \bigr)_{n\ge1}:= \Biggl(
\Biggl(\sum_{i=T(n)}^{T(n+1)-1}\tilde f(
\zeta_i),\sum_{i=T(n)}^{T(n+1)-1}a(
\zeta_i,\Delta X_i) \Biggr) \Biggr)_{n\ge1}
\]
is an $N_f$-dependent (under $\mb P$) stationary sequence with zero
means, by Lemma~\ref{l3} and the CLT for $m$-dependent sequences \cite
{Bi56}, Theorem 5.2, we conclude that as $n\to\infty$, the $\mb
P$-law of $(\frac{1}{\sqrt n}A_{\lfloor\cdot n\rfloor},\frac
{1}{\sqrt n}B_{\lfloor\cdot n\rfloor})$ converges weakly to a
Brownian motion in $\mb R^2$. Moreover, by the same argument as in
\cite{Sz2}, Theorem 4.1,
%
\begin{equation}
\label{e52} \Biggl(\lambda\sum_{i=0}^{\lfloor t/\lambda^2\rfloor}
\tilde f(\bar\zeta_i), \lambda\sum_{i=0}^{\lfloor t/\lambda
^2\rfloor}a(
\bar\zeta_i,\Delta X_i) \Biggr)_{t\ge0}
\end{equation}
converges weakly (under $\mb P$, as $\lambda\to0$) to a Brownian
motion $(\tilde{N}_t, N_t)$ in $\mb{R}^2$.
On the other hand, by (\ref{e51}) and Theorem~\ref{thm5},
%
\begin{eqnarray}
\label{e53} &&\mb E \Biggl[ \Biggl(\exp\bigl(G\bigl(t/\lambda^2,
\lambda\bigr)\bigr) \lambda\sum_{i=0}^{t/\lambda^2}
\tilde f(\bar\zeta_i) \Biggr)^{3/2} \Biggr]\nonumber
\\
&&\qquad \le\bigl(\mb E\bigl[\exp\bigl(6G\bigl(t/\lambda^2,\lambda\bigr)\bigr)\bigr]
\bigr)^{1/4} \Biggl(\mb E \Biggl[\lambda^2 \Biggl(\sum
_{i=0}^{t/\lambda^2}\tilde f(\bar\zeta_i)
\Biggr)^2 \Biggr] \Biggr)^{3/4}
\\
&&\qquad \le Ce^{ct}N_f
\llVert f\rrVert_\infty^{3/2}.
\nonumber
\end{eqnarray}
Therefore, by the invariance principle (\ref{e52}) and uniform
integrability (\ref{e53}),
%
\begin{eqnarray}
\label{e9} \lim_{\lambda\to0} \lambda\anne\Biggl[ \sum
_{i=0}^{t/\lambda^2} \bigl(f(\bar\zeta_i)-\mc Q f
\bigr) \Biggr] &=& \lim_{\lambda\to0}\mb E \Biggl[ \exp\bigl(G\bigl(t/
\lambda^2,\lambda\bigr)\bigr) \lambda\sum
_{i=0}^{t/\lambda^2}\tilde f(\bar\zeta_i) \Biggr]
\nonumber
\\
&\stackrel{\fontsize{8.36pt}{8.36pt}{(\ref{e1})}} {=}& E\bigl[\tilde N_t\exp
\bigl(N_t-EN_t^2/2\bigr)\bigr]
\\
&=&t\cov(N_1, \tilde N_1):=t\Lambda(f).\nonumber
\end{eqnarray}
Setting $U_j=\sum_{k=T(j+N_f)}^{T(j+N_f+1)-1}\tilde f(\bar\zeta_k)$ and
$V_j=\sum_{k=T(j+N_f)}^{T(j+N_f+1)-1}a(\bar\zeta_k,\Delta X_k)$,
$\Lambda$ also has the expression
%
\begin{eqnarray}
\label{e24} \Lambda(f) &=& \lim_{n\to\infty}\mb E \Biggl[\sum
_{i=0}^{T(n)}\tilde f(\bar\zeta
_i)\sum_{i=0}^{T(n)}a(\bar
\zeta_i,\Delta X_i) \Biggr] \Big/\mb E\bigl[T(n)\bigr]
\nonumber\\[-8pt]\\[-8pt]\nonumber
&=& \Biggl(\mb E[U_1V_1]+\sum
_{i=1}^{N_f}\mb E[U_1 V_{1+i}+U_{1+i}V_1]
\Biggr) \Big/\mb E\bigl[T(2)-T(1)\bigr].
\end{eqnarray}

Therefore,
\begin{eqnarray*}
&&\biggl\llvert\frac{\mc Q_\lambda f-\mc Q f}{\lambda}-\cov(N_1,\tilde N_1)
\biggr\rrvert
\\
&&\qquad \le\frac{1}{\lambda}\anne\Biggl\llvert\frac{\lambda^2}{t}\sum
_{i=0}^{t/\lambda^2} f(\bar\zeta_i)-
\mathcal{Q}_\lambda f \Biggr\rrvert+ \Biggl\llvert\frac{\lambda}{t}
\anne\Biggl[ \sum_{i=0}^{t/\lambda
^2} \bigl(f(\bar
\zeta_i) -\mc Q f \bigr) \Biggr] -\cov(N_1, \tilde
N_1)\Biggr\rrvert.
\end{eqnarray*}
Letting first $\lambda\to0$ and then $t\to\infty$, we obtain [by
Theorem~\ref{thm5} and (\ref{e9})]
%
\begin{equation}
\label{e19} \lim_{\lambda\to0}\frac{\mc Q_\lambda f-\mc Q f}{\lambda}
=\Lambda(f).
\end{equation}
Theorem~\ref{thm6} is proved.
\end{pf*}

%
\begin{remark}
1. By (\ref{e24}), for any $f$ that satisfies (\ref{f}),
\[
\bigl\llvert\Lambda(f)\bigr\rrvert\le CN_f\llVert f\rrVert
_\infty.
\]

2. By the same argument as in \cite{BZ}, one can obtain a
quenched invariance principle for (\ref{e52}).
\end{remark}

\section{The derivative of the speed and Einstein relation}\label{secE}
In this section, we will apply Theorems~\ref{ER0} and~\ref
{thm6} to derive the derivative of the speed.

\begin{pf*}{Proof of Corollary~\ref{thm1}}
Note that
\[
v_\lambda=\mc Q_\lambda\bigl[d\bigl(\omega^\lambda\bigr)
\bigr] =\mc Q_\lambda\bigl[d(\omega)\bigr]+ \lambda\mc Q_\lambda
\bigl[d(\xi)\bigr]
\]
and
\[
v_0=\mc Q\bigl[d(\omega)\bigr].
\]
Thus,
\[
\frac{v_\lambda-v_0}{\lambda} = \mc Q_\lambda\bigl[d(\xi)\bigr]+\frac
{\mc Q_\lambda[d(\omega)]-\mc Q
[d(\omega)]}{\lambda}.
\]
Therefore,
by Theorem~\ref{ER0} [recall that $\Lambda=0$ in case (i)] and
Theorem~\ref{thm6},
\[
\lim_{\lambda\to0}\frac{v_\lambda-v_0}{\lambda} =\mc Q\bigl[d(\xi)\bigr
]+\Lambda
\bigl(d(\omega) \bigr).
\]
[Here, we write $\Lambda(f):=(\Lambda f_1,\ldots, \Lambda f_d)$ for a
function $f=(f_1,\ldots,f_d)\dvtx \Omega\to\mb R^d$.]
Corollary~\ref{thm1} is proved.
\end{pf*}

\begin{pf*}{Proof of Proposition~\ref{prop7}}
The existence of the speed is proved in Proposition~\ref{prop9}. When
$\omega$ is balanced and $\xi(x,e)=\omega(x,e)e\cdot\ell$, it is
straightforward to check that $\mc Q(d(\xi))=D\ell$.
\end{pf*}

%
\begin{remark}
1. For case (ii), with Corollary~\ref{thm1}, we can also write the
derivative of the speed at $\lambda>0$:
\[
\frac{\ud v_\lambda}{\ud\lambda} = \mc Q_\lambda d(\xi)+\Lambda_\lambda
\bigl(d
\bigl(\omega^\lambda\bigr)\bigr),
\]
where $\Lambda_\lambda$ is as $\Lambda$ in (\ref{e24}), with
$\omega$, $\mb E$ and $\mc Q$ replaced by $\omega^\lambda$, $\mb
E_\lambda$ and $\mc Q_\lambda$, respectively. It is not hard (by
considering the Radon--Nikodym derivative) to obtain
\[
\lim_{\lambda\to0}\Lambda_\lambda f=\Lambda f.
\]
So $\ud v_\lambda/\ud\lambda$ is continuous at $\lambda=0$ and
hence also continuous for $\lambda\in[0,\lambda_0)$.

2. For case (i),
\[
\frac{\ud v_\lambda}{\ud\lambda}=\mc Q_\lambda d(\xi)+\lambda\Lambda
_\lambda
\bigl(d(\xi)\bigr).
\]
$\Lambda_\lambda$ can also be expressed in terms of the regeneration
times defined in Section~\ref{secreg}. Moreover, using
Lebowitz--Rost's argument and the moment estimates of the
regenerations, it is not hard to obtain $\llvert\lambda\Lambda
_\lambda(d(\xi))\rrvert\le C$. But it is not clear whether
$\lim_{\lambda\to0}\lambda\Lambda_\lambda(d(\xi))=0$, that is,
$\ud v_\lambda/\ud\lambda$ is also continuous at $\lambda=0$.

3. By (\ref{e24}), we get
\[
\Lambda(\mbox{Constant})=0.
\]
Hence when the original environment is deterministic, Corollary~\ref
{thm1} agrees with Sabot's result \cite{Sabot}. [Note that when
$\omega$ and $\xi$ are independent, $\mc Q(\xi\in\cdot)= \mc P(\xi
\in\cdot)$.]
\end{remark}
%


\section{Questions}
\begin{enumerate}
%
\item Is the Einstein relation still true for balanced environment
without the uniform ellipticity assumption?
(Recall that the quenched invariance principle for random walks in
i.i.d. balanced random environment is proved for elliptic environment
\cite{GZ} and ``genuinely $d$-dimensional'' environment \cite{BZ}.)
\item In case (i), is $\ud v_\lambda/\ud\lambda$ continuous at
$\lambda=0$? Further, is $v_\lambda$ an analytic function of~$\lambda$?
\item Does the Einstein relation hold for a random environment with
zero-speed but is not balanced, for example, RWRE with cut points \cite{BSZ}?
\item We expect Theorem~\ref{ER0} to be true for general random
environment (with an ergodic stationary measure for the environment
viewed from the particle process) with general perturbations.
But it is not clear how this can be proved.
\end{enumerate}

\begin{appendix}\label{appe}
\section*{Appendix: Proof of Theorem~\texorpdfstring{\protect\ref{44444}}{12}}
The idea of our proof is the following. Since the drift $\omega
^\lambda$ at each point is of size $c\lambda$, the ``worst case'' is
that all the drifts $d(\theta^x\omega^\lambda)$ point toward the
level $ \{z\dvtx  z\cdot e_1=0\}$. Hence, we only need to work on the
``worst case'' to get the upper bound. To this end, we couple $X_i$ with
a slow chain $Y_i$ on $\mb Z^+$, which is defined by
\begin{eqnarray*}
Y_0&=&\llvert X_0\cdot e_1\rrvert,
\\
Y_{i+1}-Y_i&=&\cases{ 0, &\quad if $X_{i+1}\cdot
e_1-X_i\cdot e_1=0$,
\vspace*{3pt}\cr
1, &\quad if
$X_{i+1}\cdot e_1\neq X_i\cdot e_1$
and $Y_i=0$,
\vspace*{3pt}\cr
B_i(X_i), &\quad if
$X_{i+1}\cdot e_1-X_i\cdot e_1=1$
and $Y_i\neq0$,
\vspace*{3pt}\cr
-1, &\quad if $X_{i+1}\cdot
e_1-X_i\cdot e_1=-1$ and $Y_i
\neq0$,}
\end{eqnarray*}
where\vspace*{1pt} $ (B_i(x) )_{i\in\mb N, x\in\mb Z^d}$ are independent
Bernoulli random variables [which are independent of $(X_j, Y_j)_{0\le
j\le i}$] such that
\[
P\bigl(B_i(x)=1\bigr)=\frac{(1-\lambda/\kappa)}{2p(x)}
\]
and
\[
P\bigl(B_i(x)=-1\bigr)=1-\frac{(1-\lambda/\kappa)}{2p(x)},
\]
where\footnote{Note that by the uniform ellipticity assumption,
\[
p(x)=\frac{\omega^\lambda(x,e_1)}{\omega^\lambda(x,e_1)+\omega
^\lambda(x,-e_1)}\ge\frac{1-\lambda/\kappa}{2}.
\]\vspace*{-10pt}}%
\[
p(x):=P^x_{\omega^\lambda}(X_1\cdot e_1=1
\mid X_1\cdot e_1\neq0).
\]
That is, $Y_i$ reflects at the origin and moves only when $X_i\cdot
e_1$ changes. When $\llvert X_i\cdot e_1\rrvert$ decreases
and $Y_i\neq0$, $Y_\cdot$ moves left. When $\llvert X_i\cdot
e_1\rrvert$ increases and $Y_i\neq0$, $Y_\cdot$ flips a coin
$B_i(X_i)$ to decide where to move.
$(Y_i)_{i\ge0}$ has the following good properties:
\begin{longlist}[2.]
\item[1.] $\llvert X_i\cdot e_1\rrvert-Y_i$ is always a
nonnegative even integer. Hence,
%
\begin{equation}
\label{e25} \tilde T_n\le\tilde S_n,
\end{equation}
where
\[
\tilde S_n:= \biggl\{i\ge0\dvtx  Y_i=\frac{n}{\lambda_1}\biggr
\}.
\]
Moreover,
%
\begin{eqnarray}
\label{e38} && P(Y_{i+1}-Y_i=\pm1\mid
X_j,Y_j,0\le j\le i)
\nonumber\\[-8pt]\\[-8pt]\nonumber
&&\qquad= \frac{1\mp\lambda/\kappa}{2}\bigl(\omega^\lambda(X_i,
e_1)+\omega^\lambda(X_i,-e_1)\bigr)
\qquad\mbox{if }Y_i>0
\end{eqnarray}
and $P(Y_{i+1}-Y_i=1\mid X_j,Y_j,0\le j\le i)=1$ if $Y_i=0$.
\item[2.]
$(Y_i)_{i\ge0}$ is not a Markov chain. But if we set $t_0=0$,
$t_{i+1}=\inf\{n>t_i\dvtx  X_n\neq X_{t_i}\}$, then
\[
Z_i:=Y_{t_i},\qquad i\ge0
\]
is a nearest-neighbor random walk on $\mb Z^+$ that satisfies
%
\begin{equation}
\label{e42} P(Z_{i+1}-Z_i=\pm1\mid Z_j, j
\le i)= \frac{1\mp\lambda/\kappa}{2} \qquad\mbox{if }Z_i>0
\end{equation}
and $P(Z_{i+1}-Z_i=1\mid Z_j, j\le i)=1$ if $Z_i=0$.
\end{longlist}

\begin{pf*}{Proof of Theorem~\ref{44444}}
Let $Y_i, Z_i, \tilde S_i, i\ge0$ be defined as above, by (\ref
{e25}), it suffices to show that for some $s>0$,
\[
E\bigl[e^{s\lambda^2\tilde S_n/n}\mid Y_0=0\bigr]<\infty.
\]
By the same argument as in \cite{Sz1}, Lemma 1.1, it is enough to show
that for any $x\in\{0,1,\ldots, n/\lambda_1\}$,
%
\begin{equation}
\label{e39} E[\tilde S_n\mid Y_0=x]\le
\frac{cn}{\lambda^2}.
\end{equation}
Putting
\[
S_n:=\inf\{i\ge0\dvtx Z_i=n/\lambda_1\},
\]
we have
\[
\tilde S_n=\sum_{i=0}^{S_n-1}t_{i+1}.
\]
Since for every $i\ge0$, $t_{i+1}$ is a geometric random variable with
success probability $\omega^\lambda(X_{t_i},e_1)+\omega^\lambda
(X_{t_i},-e_1)\ge\kappa$, we can stochastically dominate
$(t_{i})_{i\ge0}$ by a sequence of i.i.d. Geometric($\kappa$) random
variables $(G_i)_{i\ge0}$ that are independent of~$S_n$. Thus, for any
$x\in\{0,1,\ldots, n/\lambda_1\}$,
\begin{eqnarray*}
E[\tilde S_n\mid Y_0=x] &\le& E[S_n/\kappa
\mid Z_0=x]
\\
&\le& E[S_n\mid Z_0=0]/\kappa.
\end{eqnarray*}
Therefore, to prove (\ref{e39}), we only need to show that
%
\begin{equation}
\label{e41} E[S_n\mid Z_0=0]\le\frac{cn}{\lambda^2}.
\end{equation}
With abuse of notation, we write $P(\cdot\mid Z_0=x)$ and $E[\cdot
\mid Z_0=x]$ as $P^x(\cdot)$ and $E^x[\cdot]$, respectively.

Set
\[
H_n=\inf\bigl\{i>0\dvtx  Z_i\in\{0,n/\lambda_1
\}\bigr\}.
\]
Conditioning on the hitting time to the origin, we have
%
\begin{equation}
\label{e40} E^0[S_n]=1+E^1[H_n]+P^0(Z_{H_n}=0)E^0[S_n].
\end{equation}
By (\ref{e42}), $Z_m-Z_0-m\lambda/\kappa$ is a martingale for $0\le
m\le H_n$. Thus, by the optional stopping theorem, for any $x\in\{
1,\ldots,n/\lambda_1\}$,
\[
E^x \biggl[Z_{H_n}-x-\frac{\lambda}{\kappa}H_n
\biggr]= 0.
\]
Hence,
%
\begin{equation}
\label{e43} E^1[H_n]\le\kappa E^1[Z_{H_n}]/
\lambda.
\end{equation}
By (\ref{e42}) and Proposition~\ref{prop6}, we get
%
\begin{equation}
\label{e44} c\lambda\le P^1(Z_{H_n}=n/
\lambda_1) \le c' \lambda
\end{equation}
and so
\[
E^1[Z_{H_n}]\le c\lambda\cdot n/\lambda_1\le
cn.
\]
This and (\ref{e43}) yield
\[
E^1[H_n]\le cE^1[Z_{H_n}]/\lambda
\le cn/\lambda.
\]
It then follows by~(\ref{e40}) and (\ref{e44}) that
\[
E^0 [S_n] =\frac{1+E^1[H_n]}{P^0(Z_{H_n}=n/\lambda_1)} = \frac
{1+E^1[H_n]}{P^1(Z_{H_n}=n/\lambda_1)} \le
C\frac{n}{\lambda^2}.
\]
Inequality (\ref{e41}) is proved.
\end{pf*}
\end{appendix}

%
\section*{Acknowledgments}
I thank my PhD advisor, Ofer Zeitouni, for introducing this problem to
me and for numerous useful discussions.
I thank Noam Berger, Nina Gantert and Pierre Mathieu for many
interesting discussions on the Einstein relation.


%

\printaddresses
\end{document}